\documentclass[12pt]{article}

\usepackage{amsmath,amsthm,amssymb}

\theoremstyle{plain}
\newtheorem{Thm}{Theorem}
\newtheorem{Prop}{Proposition}
\newtheorem{Cor}{Corollary}
\newtheorem{lem}{Lemma}[section]

\makeatletter

\@addtoreset{equation}{section}

\def\R{\mathbb{R}}
\def\Z{\mathbb{Z}}

\def\v2{\vskip2mm}
\def\n{\noindent}

\def\({(\!(}
\def\){)\!)}
\def\a{\alpha}
\def\b{\beta}
\def\e{\varepsilon}
\def\de{\delta}

\def\la{\lambda}

\def\Om{\Omega}

\def\pf{{\it Proof.}}
\def\v2{\vskip2mm}
\def\n{\noindent}

\def\0{{\bf 0}}

\def\tst12{{\textstyle \frac12}}

\def\n{\noindent}
\def\beq{\begin{eqnarray*}}
\def\eeq{\end{eqnarray*}}

\def\beqn{\begin{equation}}
\def\eeqn{\end{equation}}

\begin{document}

\begin{center}
{\bf One dimensional  random walks  killed on a finite set} \\
\vskip4mm
{K\^ohei UCHIYAMA} \\
\vskip2mm
{Department of Mathematics, Tokyo Institute of Technology} \\
{Oh-okayama, Meguro Tokyo 152-8551\\
e-mail: \,uchiyama@math.titech.ac.jp}
\end{center}

\vskip8mm
\n
{\it running head}:   random walk  killed on a finite set

\vskip2mm
\n
{\it key words}:  
  exterior domain; transition probability; escape from a finite set; hitting probability of a finite set;  potential theory
\vskip2mm

{\it AMS Subject classification (2010)}: Primary 60G50,  Secondary 60J45. 

\vskip6mm

\begin{abstract}
We study the   transition probability,  say  $p_A^n(x,y)$,   of a  one-dimensional random walk on the integer lattice killed when entering into a non-empty finite set $A$.  The random walk is assumed to be  irreducible and have zero mean and a finite variance $\sigma^2$. We show that  $p_A^n(x, y)$ behaves  like 
$[g_A^{+}(x)\widehat g_{A}^{\,+}(y) + g_A^-(x)\widehat g_{A}^{\,-}(y)] (\sigma^{2}/2n) p^n(y-x)$ uniformly  in the regime characterized by  the conditions $|x|\vee |y| =O(\sqrt n)$ and  $|x|\wedge |y|= o(\sqrt n)$ generally  if $xy>0$ and under a mild additional assumption about the walk if $xy<0$.  Here $p^n(y-x)$ is the transition kernel of the random walk  (without killing);  $g^\pm_A$ are  the Green functions for the \lq exterior' of $A$ with  \lq pole at $\pm \infty$'  normalized so that $g^\pm_A(x) \sim 2|x|/\sigma^2$ as $x \to \pm\infty$;  and $\widehat g_A^{\, \pm}$ are the corresponding Green functions for the time-reversed  walk. 
\end{abstract}
\vskip6mm

\begin{center} {\sc Contents}
\end{center}

1.\,  Introduction and main results

2.\, Preliminaries from the theory of one-dimensional random walks

3.\,  The Green function and escape from A

\quad 3.1.\, Functions $u_A$, $g_A^+$ and $g_A^-$

\quad 3.2.\, Probabilities of escape from $A$ and an overshoot estimate

4.\, Results for a single point set and a half line

\quad 4.1.\, Results for the case $A=\{0\}$

\quad 4.2.\, Space-time distribution of entrance into $(-\infty,0]$

5.\, Proof of  Theorem 1

6.\, Refinements in  case $xy<0$

7.\, Appendices

\quad A. \, A consequence of duality

\quad B. \, Some inequalities concerning $ a(x) -x/\sigma^2$

\quad C. \, Comparison between $p_A^n$ and  $p_{\{0\}}^n$

References

\vskip6mm

\section { Introduction and main results}

This paper concerns the transition probability of  a one-dimensional  random walk on the integer lattice  $\Z$ killed on a finite set $A$. For  random walks on the $d$-dimensional integer  lattice $\Z^d$, $d\geq1$ killed on a finite set 
 H. Kesten  \cite{K} obtained,   among others,  the  asymptotic form of the transition probability under a quite general setting.  For the  important case of one dimensional random walk with zero mean and  finite variance, however, his result  is restricted to the special case when  $A$ consists of a single  point. In this paper 
 we extend it to every finite set $A$.  Our result is stronger than his in another respect: the asymptotic estimate is valid uniformly for space variables within a reasonable  range of relevant variables.   It is incidentally  revealed  that according as the third absolute moment of the increment variable is finite or not,  the walk killed at the origin  exhibits qualitatively different behaviour as the starting and landing 
 positions are taken far from the origin in the opposite  directions from each other (see Remark 2 near the end of this section).
 
Our method of proof  is quite different from that of Kesten \cite{K}. In \cite{K} a compactness argument is 
used as a basic tool. Our proof  reflects the behaviour  of random walk path. It   rests on the results of \cite{U1dm} in which  the same problem as the present paper is studied but with $A=\{0\}$ and an asymptotic form of the transition probability  valid uniformly  for  space variables is  obtained.  
The same method is applied in  \cite{Uf.s} to   higher dimensional random walks to obtain  a similar strengthening of Kesten's result.  For multidimensional Brownian motions the corresponding problem is studied by  \cite{CMM} for space variables restricted to compact sets and by \cite{Ubdsty} without  any restriction as such.

Let $S_n=S_0+X_1+\cdots+X_n$, $n=1,2,\ldots$ be a  random walk on the one-dimensional integer lattice  $\Z$. Here the  increments $X_j$ are i.i.d.  $\Z$-valued random variables  defined on some probability  space $(\Om, {\cal F}, P)$ and the initial state $S_0$ is an integer   left unspecified for now. 
As usual the  law with  $S_0=x$  of the walk  $(S_n)$  is denoted by  $P_x$ and the corresponding expectation by $E_x$.   Throughout this paper we suppose that the random walk  $(S_n)$ is irreducible, namely  for every $x\in \Z$,     $P_0[ S_n=x]>0$ for some $n>0$, and    that  
\beqn\label{mom}
EX=0~~~~\mbox{and}~~~~0<\sigma^2 := E X^{2}<\infty  
\eeqn
Here as well as in what follows  $X$ is a random variable having the same law as $X_1$ and $E$ 
the expectation w.r.t. $P$.

Let $p(x) = P[X=x]$ and   $p^n(x)= P_0[S_n=x]$ so that for $y\in \Z$,
\beq
P_x[S_n=y] =p^n(y-x) \quad\mbox{and} \quad p^0(x)=\de_{x,0},
\eeq  
 where $\de_{x,y}$ equals unity if $x=y$  and zero  if $x\neq y$. 
For a  non-empty finite subset $A$ of $\Z$, let $p_A^n(x,y)$ denote the transition probability of the walk $S_n$ killed upon entering $A$, defined by
$$p_A^n(x,y) = P_x[ S_k\notin A\,\,\mbox{for}\,\,1\leq  k\leq  n\,\, \mbox{and}\,\, S_n =y],  \quad n=0,1, 2,\ldots.$$
Thus $p_A^0(x,y)= \de_{x,y}$ (even if $y\in A$); and $p_A^n(x,y) =0$ whenever $y\in A, n\geq 1$.

Let $a(x)$ be 
the potential function  of the walk defined by
$$a(x) = \lim_{n\to\infty}\sum_{k=0}^n [p^k(0) - p^k(-x)].$$
It is convenient to bring in 
$$a^\dagger(x) :=\de_{x,0} + a(x).$$
The result of Kesten \cite{K} mentioned above implies   that for each  $x$ and $y\neq 0 $, as $n\to\infty$ 
  \beqn\label{K}
   \lim_{n\to\infty} \frac{p_{\{0\}}^n(x,y)}{f_0(n)} = a^\dagger(x)a^\dagger(-y)+\frac{xy}{\sigma^4},
   \eeqn
provided that  the walk is (temporally) aperiodic in addition.   Here 
 $$f_{0}(n) = P_0[ S_k\neq 0\,\,\mbox{for}\,\, k=1,\ldots, n-1\,\, \mbox{and}\,\, S_n =0].$$
We know that
  $$f_0(n) = \frac{\sigma^2}{n}p^n(0)\{1+o(1)\}$$
as $n\to \infty$.   
(Cf. \cite{S}  for the existence of $a(x)$ and \cite[Theorem 8]{K}  for the asymptotic form of $f_0(n)$ stated above.) Note that $y=0$ is reasonably excluded in (\ref{K}) (cf. Remark 1 (e)).

Denote  the Green function of the killed walk 
 by $g_A(x,y)$: 
$$g_A(x,y) = \sum_{n=0}^\infty p_A^n(x,y),\quad  x, y\in \Z.$$
(To be precise this  does not conform  to the usual nomenclature,  according  to which a Green function is set   zero on $A\times A$, while our   $g_A(x,y)$ is equal to $ \de_{x,y}$ if $y\in A$ and  to  
the  probability that  the first entrance    into  $A$ takes place at $x$ for  the dual walk starting at $y$ if $x\in A, y\notin A$.)
According to  Theorem 30.1 of \cite{S} $g_A(x,y)$ has limits as $y\to +\infty$ and $y\to -\infty$ for each  $x$. We call them 
$$g^+_A(x)= \lim_{y\to\infty} g_A(x, y) \quad \mbox{and} \quad g^-_A(x)= \lim_{y\to -\infty} g_A(x, y).$$
($g^\pm_A(x)$ may be interpreted as the expected number  of visits  to $x$ made by  the dual\textemdash or time-reversed\textemdash  random walk \lq started at $\pm \infty$' up to and including  the first entrance 
time into   $A$.) 
 We write  $-A$ for  $\{-z: z\in A\}$; $s\wedge t$ and $s\vee t$ for the minimum and maximum, respectively,  of real numbers $s$ and $t$.  In the following theorem we impose on the pair of  $x, y$ the condition
 \beqn\label{nat_c}
 g_A^+(x)g_{-_A}^-(-y) + g_A^-(x)g_{-A}^+(-y) \neq 0.
 \eeqn

\begin{Thm}\label{thm1} Let  $A$ be a  non-empty finite subset of $\Z$. 

{\rm (i)}  \, For each $M\geq 1$,  uniformly for $x\in \Z$ and $y\in \Z\setminus A$ subject to condition (\ref{nat_c}) and the constraints  $-M\leq x\leq M\sqrt  n$  and $-M\leq  y\leq M\sqrt  n$, as  $n\to\infty$ and $(|x|\wedge |y|)/\sqrt n \to 0$
\beqn\label{eq_thm1}
p_A^n(x,y) = \frac{g_A^+(x)g_{-_A}^-(-y) + g_A^-(x)g_{-A}^+(-y)}{2n/\sigma^2}p^n(y- x)\{1+o(1)\}; 
\eeqn
if condition (\ref{nat_c}) is violated, then $p_A^n(x,y) < C e^{-\la n}$ for some positive constants $\la$ and $C$ that depend only on $p$ and $A$.  

{\rm (ii)}  \,  As $x\wedge y \wedge n \to\infty$ under $x\vee y<M\sqrt n$ along with  $p^n(y-x)>0$
\beqn\label{eq_thm10}
p_A^n(x,y) = \frac{\nu}{\sqrt{2\pi \sigma^2 n}}\Big(e^{- (y-x)^2/2\sigma^2 n} - e^{- (y+x)^2/2\sigma^2 n}\Big)\{1+o(1)\},
\eeqn
 where $\nu$ designates the (temporal) period of the random walk. 
\end{Thm}

\v2
In the theorem above as well as  in what follows,  $o(1) \to 0$ in the specified  procedure of taking limit and the convergence is uniform under the specified constraint.  By symmetry the results of Theorem \ref{thm1}
are valid if $x$ and $y$ are simultaneously replaced by $-x$ and $-y$, respectively, and  this  remark applies to the succeeding results.   The case $xy<0$ and $|x|\wedge |y|\to\infty$,   excluded in  Theorem \ref{thm1},
is discussed  later in this introduction;  as a matter of fact  formula (\ref{eq_thm1}) remains true under a mild additional assumption about $p$  but may break down without it (cf. Theorem \ref{thm2} below and  Theorem \ref{thm4} and  Remark 6 in Section 6).  
    \v2
{\sc Remark 1.} (a) \, If the object corresponding to the dual (time-reversed) walk is indicated by putting $\, \widehat\,$\,,  like $\widehat p^{\,\, n}_A(y,x)$,  we have for $x, y\notin A$
$$ g_{-A}(-y,-x) = \widehat g_A(y,x) = g_A(x,y), $$
hence
 \beqn\label{dual1}
 g_{-A}^-(-y) =\widehat g_A^{\,+}(y) =\lim_{x\to\infty}g_A(x,y)
 \eeqn
 (cf. \cite[Section 10]{S}; see also  Appendix A for explanation and related matters). 
 \v2
 
(b)\, We shall show  
$$g^\pm_A(x) = a(x)\pm x/\sigma^{2} + O(1),$$
 (see (\ref{repr_g}), (\ref{repr_g-})). Here as elsewhere both upper or both  lower signs should be chosen in the 
 double signs. 
It in particular follows that 
$$g^+_A(x)/x \, \longrightarrow\, 2/\sigma^2 \,\,\, \mbox{or}\,\,\, 0 \quad\mbox{according as}
\quad x\to +\infty\,\,\,  \mbox{or}\,\, -\infty;$$
and similarly for $g^-_A(x)$ and  $\widehat g^{\,\pm}_A(y)$.  

By substitution of these relations  formula (\ref{eq_thm1}) is somewhat simplified when $x\vee y \to \infty$ (under $x\vee y <M\sqrt n$ and  $x\wedge y = o(\sqrt n)$). 
In fact, in the case  $y\to\infty$, we have $g^+_{-A}(-y) =o(y)$ as well as $g^-_{-A}(-y)\sim y/\sigma^2$ and,  if (\ref{nat_c}) holds,    $p_A^k(x,y)>0$ for some $k\geq 1$, which implies  $g^+_A(x)>0$ (cf. (\ref{V})), so that (\ref{eq_thm1}) can be written as 
 \beqn\label{R1(c)}
 p_A^n(x,y) \sim \frac{g^+_A(x)y}{n}p^n(y- x),
 \eeqn
 where the symbol  $\sim$ means that the ratio of two sides of it approaches unity;
and similarly for the case $x\to\infty$ (cf. Lemma \ref{lem4.2}).

Also, as $ x\wedge y\to\infty$ (under $x\vee y <M\sqrt n$ and $x\wedge y=o(\sqrt n\,)$), (\ref{eq_thm1})  reduces to
$$
p_A^n(x,y) \sim \frac{2xy}{\sigma^2n}p^n(y- x).
$$
 This relation  conforms to  (\ref{eq_thm1}) in view of a local central limit theorem. It also follows that 
 the restriction $x\wedge y = o(\sqrt n)$ cannot be relaxed in (\ref{eq_thm1}). 
Generally  in the parabolic regime $|x|\vee |y| =O(\sqrt n)$ it holds  that 
 $$p_A(x,y) \asymp |xy|n^{-1} p^n(y-x)  \quad\mbox{if}   \quad xy>0$$
($\asymp$  means that the ratio of two sides is bounded away from zero  and infinity) and 
  $$C( |x|+ |y|) n^{-1} p^n(y-x) \leq p_A(x,y) = p^n(y-x) \times o(|xy|/n) \quad\mbox{if}   \quad xy<0,$$
 provided  condition (\ref{nat_c}) is satisfied. (See  Theorem \ref{thm2} and Remark 2 after it for the latter.)
  \v2
  
   (c) \,  Because  of (\ref{dual1}) $g^-_{-A}(-y)$ (resp. $g^+_{-A}(-y)$) is positive  if and only if $y$ can be reached by the walk starting at  $+\infty$ (resp. $-\infty$). Taking  account of  this and its analogue for $g^\pm_A$ we see that  the product  $g^+_A(x)g^-_{-A}(-y)$ (resp. $g^+_A(x)g^-_{-A}(-y)$) is positive  if (and only if)
 the walk starting at  $x$ can  reach  $y$ after a large excursion to the right (resp. left) with a positive probability.  Hence,  if both of  these two products vanish, namely  if condition (\ref{nat_c})  is violated, then one of the following (1) or (2) must hold true:
 \v2
(c1) \;\; $p_A^n(x,y)$ vanish for all  $n$;  
 
(c2) \;\;  all the random walk paths up to $\sigma_A$ that start at $x$ and pass through  $y$  are  \\
\qquad\qquad\;\,  confined in the convex hull of $A$ with probability one,
\v2 \n
where $\sigma_A$ denotes the first entrance time into $A$  (see (\ref{sigm/tau}) below for the  precise definition).
 In the second case the paths  must enter $A$ in a small number of steps so as to ensure  the second 
 assertion of Theorem \ref{thm1}.

  \v2
(d)\,  Formula (\ref{eq_thm1}) may be equivalently stated as follows:  The probability that the walk $(S_n)$ started at  $x$ and pinned at $y$ at time  $n$ avoids  $A$ is asymptotically equivalent  to the ratio on the right side of (\ref{eq_thm1}), namely
$$P_x[\sigma_A >n\,|\, S_n=y] \, \sim  \, \frac{g_A^+(x)g_{-A}^-(-y) + g_A^-(x)g_{-A}^+(-y)}{2n/\sigma^2}$$
for $n, x, y$ such that $p^n(y-x)>0$. 

 \v2
 (e)\, Formula (\ref{eq_thm1}) implies Kesten's result (\ref{K}). Indeed, $g_{\{0\}}^\pm(x)= a^\dagger(x)\pm x/\sigma^2$ so that for $A=\{0\}$ the numerator of the ratio on the right side of (\ref{eq_thm1})  reduces to $a^\dagger(x)a^\dagger(-y) +xy/\sigma^4$ (cf. (\ref{1*}), (\ref{repr_g})).
What is actually treated, instead of our $p_A^n(x,y)$, in \cite{K} is
$$Q_A^n(x,y) := P_x[ S_k\notin A\,\,\mbox{for}\,\, 0< k< n\,\, \mbox{and}\,\, S_n =y],  \quad n=0,1, 2,\ldots,$$
so that $Q_A^n(x,y)$ equals $ p_A^n(x,y)$ or $P_{x}[ \sigma_A=n, S_{\sigma(A)}= y]$ according as $y\notin A$ or $y\in A$. If  $p^n_{\{0\}}(x,y)$ is replaced by $Q^n_{\{0\}}(x,y)$, then formula (\ref{K})  
 becomes valid also for $y=0$. The same remark applies to  (\ref{eq_thm1})  for  $y\in A$ as formulated below in Corollary \ref{cor1}. The corresponding Green function $G_A(x,y) :=\sum_{n=0}^\infty Q_A^n(x,y)$ is related to $g_A$ by  
 \beqn\label{G-g}
  G_A(x,y) =g_A(x,y)+ P_{x}[S_{\sigma(A)}= y].
  \eeqn
  Although $Q_A^n$ and  $G_A$ are natural objects to consider because of their  symmetry relative to duality, we adhere to $p^n_A(x,y)$ as the  principal object to study.

 \v2
 (f)\,    For random walks with drift (i.e. in case $\mu:=E X \neq 0$)  one can readily derive from Theorem \ref{thm1}   the corresponding asymptotic formulae (namely those in the regime $|x|\vee |y|=O(\sqrt n$) if $\sum p(x)s_0^xx=0$ and  $\sum p(x)s_0^xx^2 <\infty$ for some $s_0\neq 1$.  In the case when there exists no such $s_0$ Kesten's results \cite[Theorems 5 and 6]{K} provide a certain   asymptotic formula under an additional  condition. 
  From another view point it is natural to take up the regime  $ |y-x-\mu n| \leq M\sqrt n$. Without difficulty   one can show that  if  $\mu> 0$,  then for each $\de>0$, uniformly  for $x \geq (- \mu +\de) n$  as $n\to\infty$ and  $y$ in this regime
 \[\label{Eq.P70} 
 p_A^n(x,y) \sim P_x[\sigma_A=\infty] p^n(y-x).
 \]
  Some detailed  investigation  is  undertaken  in a separate paper \cite{U1dmk}.
 
 \v2
 (g)\,  Extensions to  non-lattice walks would be interesting. The present work rests on the results for the case $A=\{0\}$ given in  \cite{U1dm}  for which  the harmonic analysis  is effectively applicable but  does not seem to be for the non-lattice walks.  However, another approach seems promising to lead to  corresponding results at least under some assumption on the distribution of the increment variable:  one may take the half line $(-\infty,0]$ in place of $\{0\}$ (see Remark 4 of Section 5)  and apply the  corresponding results for the transition probability that are found in  several papers \cite{BD}, \cite{Car}, \cite{D}, \cite{VW}, etc. and those for the  potential operator given in \cite{PS}.
 
 \v2\v2

For a non-empty set $B\subset \Z^d$,  $\sigma_B$ (resp. $\tau_B$)  denotes the first time when  $S_n$  enters into (resp. exits from) $B$:
\beqn\label{sigm/tau}
  \sigma_B = \inf\{n\geq 1: S_n \in B\},\quad\quad  \tau_B = \inf\{n\geq 1: S_n \notin B\}. 
  \eeqn
For typographical reason  we shall sometimes write $\sigma(B)$ for $\sigma_B$ and similarly for $\tau(B)$. 
 For a  positive integer $R$ denote the interval  $\{-R+1, \ldots, R-2,  R-1\}$ by 
 $$U(R).$$ 

 The roles of $g^{\pm}_A$ and $g^{\pm}_{-A}$ appearing in (\ref{eq_thm1}) may be explained by  the formula
\beqn\label{crucial}P_x[\tau_{U(R)} =\sigma_{[R,\infty)} <\sigma_A] = \frac {g^+_A(x)}{R}\{1+o(1)\}
\eeqn
  (Proposition \ref{prop3} (\S2))   and its reverse and dual formulae  (see the  paragraph  given near the end of this section for a little more details). Here $o(1) \to 0$  as  $R\to\infty$ uniformly for $-M < x \leq R$. Indeed formula  ({\ref{crucial})  is a gist of the proof of Theorem \ref{thm1},  and   Section 3 will be devoted to its proof. 
\v2

For $\xi_0\in A$ and $x\in \Z$, put
$$w_A(x) = \sigma^{-2}( x - E_x[S_{\sigma(A)}])$$
and 
$$u_A(x) = a^\dagger(x-\xi_0) -E_x[a(S_{\sigma(A)} - \xi_0)], $$
where the right side does not depend on the choice of $\xi_0$ (cf. \cite{Uf.s}, see also (\ref{2_1})).
  We shall see  that
\beqn\label{1*}
g^+_A(x)= u_A(x) +w_A(x)\quad\mbox{and} \quad g^-_A(x)= u_A(x) -w_A(x),
\eeqn
by which we obtain the identity
$$\frac{g^+_A(x)g^-_{-A}(-y) + g^-_A(x)g^+_{-A}(-y)}2= u_A(x)u_{-A}(-y) - w_A(x)w_{-A}(-y).$$
By substitution  formula (\ref{eq_thm1})  becomes quite analogous to  the one for the case $A=\{0\}$ as given in Theorem A (i) of Section 4.1 that  extends   (\ref{K}) to  unbounded $x,y$.

Noting that $g_A(\cdot, y)$ is bounded for each  $y$,  we  pass to the limit in the identity
$$\sum_{z\in \Z} p^1_A(x,z)g_A(z,y) =g_A(x,y) - \de_{x,y}$$
to find  that $g^+_A$ is a harmonic function for the killed walk in the sense that
\beqn\label{harm_g}
g^\pm_A(x) = \sum_{z\in \Z\setminus A} p(z-x) g^\pm_A(z)   \quad\mbox{for all}\quad  x\in \Z;
\eeqn
the functions  $u_A$ and $w_A$ also are harmonic in the same sense, and $g_{-A}^\pm(-\,\cdot)$, $u_{-A}(-\,\cdot)$ and $w_{-A}(-\, \cdot)$ are dual harmonic.

 From the proof of Theorem \ref{thm1} or directly by a usual argument based on the last leaving decomposition with the help of  the  relations dual to  (\ref{harm_g}) (see Remark 5 at the end of Section 5)  we obtain  the following  
\begin{Cor}\label{cor1} Let  $A$ be a  non-empty finite subset of $\Z$. Then,  for each $M>1$ and  for $\xi\in A$, 
as $n\to\infty$ 
\beqn\label{2}
P_x[ \sigma_A =n, S_{\sigma(A)} =\xi] = \frac{g^+_A(x) g^-_{-A}(-\xi) +g^-_A(x) g^+_{-A}(-\xi)     }{2n/\sigma^2}p^n(\xi-x) \{1+o(1)\}
 \eeqn 
  uniformly for  $x\in \Z$ satisfying 
$|x|< M\sqrt n$ and  $g^+_A(x) g^-_{-A}(-\xi) +g^-_A(x) g^+_{-A}(-\xi) >0$; if the latter condition is
 violated,  then the probability on the left side of (\ref{2}) either vanishes for every $n$ or tends to zero exponentially fast as $n\to\infty$.
\end{Cor}

\v2
For a non-empty  set $B$ that is contained in $(-\infty, N]$ for some  $N$, we put 
 $$H^{+\infty} _B (y) = \lim_{x\to\infty} P_x[S_{\sigma(B)} =y]$$ 
and similarly for $H^{-\infty}_B$ if  $B \subset  [N,\infty)$ (the limits exist and  $H^{\pm \infty}_B$ are probabilities as  is established 
in \cite[Theorem 30.1]{S}).
It is noted that  $g^+_{A}(\xi)=\sum_{z\notin A} p(z -\xi)g^+_{A}(z)$ is the probability that the dual walk \lq started at $+\infty$' hits $A$ at $\xi$, which fact  is expressed as
$$g^+_{A}(\xi)=  H_{-A}^{-\infty}(-\xi), \quad \xi\in A,$$
and similarly   $g^-_A(\xi) = H^{+\infty}_{-A}(-\xi)$, representing the same probability but  for the dual walk started at $-\infty$.  Thus  $g^{-}_{-A}(-\xi)$ is the probability that the  dual walk  \lq started at $- \infty$' hits $-A$ at $-\xi$,  whence
\beqn\label{g_pm}
 g^{-}_{-A}(-\xi) = H^{+\infty}_A(\xi) ,\quad \xi\in A,
\eeqn
and similarly   $g^{+}_{-A}(-\xi) = H^{-\infty}_A(\xi)$.
By (\ref{g_pm})  and by (\ref{1*})  we especially  have
$$\sum_{\xi\in A}g_A^{+}(\xi) =  \sum_{\xi\in A}g_A^{-}(\xi) =1. 
$$
It is noted that 
$$ \frac{g^+_A(x)+g^{-}_A(x) }2=u_A(x).
$$

The asymptotic form of $P_x[\sigma_A=n]$ is obtained by making summation over $\xi\in A$ in  (\ref{2}). In general, however, we need to take care  of the temporal periodicity of the walk and partition the lattice  $\Z$  according to it. In the following corollary we suppose for simplicity  that 
the random walk is  aperiodic (strongly  aperiodic in the sense of \cite{S}) so that $p^n(\xi-x)$ may be replaced by $p^n(-x)$. 

\begin{Cor}\label{cor12}  Suppose  the random walk is  aperiodic. Then 
$$ P_x[ \sigma_A =n] = \frac{\sigma^2 u_A(x)}{n}p^n(-x) \{1+o(1)\}$$
as $n\to \infty$; and for $\xi\in A$,
$$\lim_{n\to \infty} P_x[S_{\sigma(A)}=\xi \,|\, \sigma_A =n] 
=  \frac{g^+_A(x)}{g^+_A(x) + g^-_A(x)} H^{+\infty}_{A}(\xi) +\frac{g^-_A(x)}{g^+_A(x)+g^-_A(x)} H^{-\infty}_{A}(\xi).$$
Here $o(1)$ and the convergence in the $\lim_{n\to \infty}$ are uniform for $|x|< M\sqrt n$.
\end{Cor}

\v2
It is often useful to have an upper bound of $p^n_A(x,y)$ valid for all $x,y$. The following one is essentially a corollary of  \cite[Theorem 1.1]{U1dm}  (see the end of Section 4.1 for the proof).
\begin{Prop}\label{prop1.1} There exists a constant $C$ such that if $|x|\wedge|y|\wedge n\geq 1$,
\beqn\label{1.8} 
p_A^n(x,y) \leq \frac{|x|\vee|y|}{|x|\vee|y|\vee \sqrt n}\Bigg[C\frac{(|x|\wedge|y|) {\sf g}_{4n}(y-x)}{(|x|\vee|y|\vee \sqrt n)} + o\bigg(\frac{|x|\wedge |y|\wedge \sqrt n\,}{(y-x)^2\vee n}\Bigg)\bigg],
\eeqn
where $ {\sf g}_n(t)= (2\pi\sigma^2 n)^{-1/2} e^{-t^2/2\sigma^2n}$.
\end{Prop}

\v2
The case when  $y<0<x$ and  $x\wedge |y| \to \infty$---excluded from Theorem \ref{thm1}---requires  somewhat delicate analysis to find whether  formula (\ref{eq_thm1}) holds true.  We shall prove in Theorem  \ref{thm4} 
that  (\ref{eq_thm1})  remains valid  for $x>0, y<0$ satisfying the  constraints $|x|\vee|y|<M\sqrt n$ and $|x|\wedge|y|=o(\sqrt n\,)$ under a mild additional condition on $p$ which is true at least if the tail $F(y)= P[X<y] $ is regularly varying at $-\infty$ with an exponent less than $-2$.  For its proof 
we need to know a certain property of the difference   $a(z)-z/\sigma^2$ (not $a(z)+z/\sigma^2$ as in the case $y>0$) for large positive values of $z$, which is  sensitive to the behaviour of the tail $P[ X< -x]$ for large $x$ in case  
$E[|X|^3; X<0] =\infty$ (cf. Lemma \ref{lem3.1}).  
Under the assumption $E[\, |X|^3; X<0] <\infty$, however,  things are particularly simplified  as given below. 
 
 Put
 $$C^+: =  \sum_{y=-\infty}^0 H^{+\infty}_{(-\infty,0]}(y) (\sigma^2a(y) +|y|) \leq \infty.$$
  It holds that
 \beqn\label{C_+}
 C^+ = \lim_{x\to +\infty} (\sigma^2a(x)-x),
 \eeqn
 $C^+<\infty$ if and only if $E[\, |X|^3; X<0] <\infty$,  and  $C^+>0$ unless 
 the walk is {\it left-continuous}, i.e., unless $P[X\leq -2] =0$ (cf. Corollary 2.1 and Remark 1 (a) of \cite{U1dm}).
 We call 
 \beqn\label{C+} 
 C^+_A := \sigma^2  \lim_{x\to +\infty} g^-_A(x) = C^+ - \sum_{\xi\in A} H^{+\infty}_A(\xi)[\sigma^{2} a(\xi-\xi_0) - (\xi-\xi_0)].
 \eeqn
 (Of course  $\xi_0\in A$; the sum does not depend on the choice of $\xi_0$ (see (\ref{D++})).) It will be proved that  $C^+_{-A} =C^+_{A}$, namely
  \beqn\label{C-}
  C^+_A =\sigma^2 \lim_{x\to +\infty} g^-_{-A}(x)  
  \eeqn
  (Lemma {\ref{lem7.1}) and that $C_A^+>0$ unless $g^-_A(x)$ vanishes for all sufficiently large (positive) $x$ (Lemma \ref{lem2.3}).   By (\ref{1*}) $\sigma^2 g^+_A(x)\sim 2x$ ($x\to\infty$) and  $\sigma^2 g^+_{-A}(-y)\sim 2|y|$ ($y\to -\infty$). Putting these together yields that  if $C_A^+<\infty$, then as $x\wedge (-y)\to\infty$
 $$\frac{g_A^+(x)g_{-A}^-(-y) + g_A^-(x)g_{-A}^+(-y) }2=C_A^+\frac{ x+ |y|}{\sigma^4}\{1+o(1)\}. $$
  In view of this relation the next result is a natural extension of Theorem \ref{thm1} to the case $(-y)\wedge x\to\infty$;  it  also extends  Theorem 1.2 of \cite{U1dm} where the same result is obtained for  $A=\{0\}$.
\begin{Thm}\label{thm2}
 If  $E[ \,|X|^{3}; X<0]<\infty$, then  as $x\wedge (-y) \wedge n\to\infty$ subject to the condition  $x\vee(- y) <M\sqrt n$
$$p_A^n(x,y) =C_A^+ \frac{x+ |y|}{\sigma^2n}p^n(y-x)\{1+o(1)\}.$$
\end{Thm}
\v2

\v2
{\sc Remark 2.}  
(a)\; In Section 6 we shall see  the following results: within the parabolic regime $-M\sqrt n < y<0 < x < M\sqrt n$
\v2
(a1)\;\;  $p_A^n(x,y)  \leq C\{g_A^+(x)g_{-A}^-(-y) + g_A^-(x)g_{-A}^+(-y) \}n^{-3/2}$
\v2
(a2)\;\;  $p_{\{0\}}^n(x,y)\geq C\{ (x+ |y|) n^{-3/2} \} \sum_{w=2}^{x\wedge |y|}p(-w) w^3$\quad ($x\wedge |y|\geq 2$)
\v2\n
(according to Propositions \ref{prop60} and \ref{prop6}, respectively).
\v2

(b)\; Suppose that $0\in A$ and  $\sharp A\geq 2$ ($\sharp A$ denotes the cardinality of $A$)  and ask in what situations  $p_A^n(x,y)$ exhibits  different asymptotic behaviour  than
$p_{\{0\}}^n(x,y)$ does. Suppose  the walk is aperiodic for simplicity and restrict the variables to the regime $|x|\vee |y|=O(\sqrt n)$. 

First   note that $p_A^n(x,y) = p_{\{0\}}^n(x,y)$ ($y\neq 0$)  for all $n$ if and only if every
path of the walk from $x$ to $y\neq 0$ that avoids $0$ also avoids $A$ ($y=0$ is reasonably excluded).   Preclude from our consideration  this   trivial situation and the case when $p^n_{\{0\}}(x,y) =0$ for all $n\geq 1$. Then,  as $|x|\vee |y| \vee n\to \infty $  under the condition  $|x|\wedge |y| =O(1)$   the ratio  $p_A^n(x,y) / p_{\{0\}}^n(x,y)$ converges to  a function of $x$ or $ y$ or  of both that is less than 1, 
as one may expect and deduces from   Theorem \ref{thm1}   (see Appendix C for a proof and more details). On the other hand   if $xy >0$ and  $|x|\wedge |y| \wedge n\to\infty$, $p_A^n(x,y)\sim p_{\{0\}}^n(x,y)$, so that the asymptotic form of $p_A^n(x,y)$  does not depend on $A$ in this regime  as is obvious from Theorem \ref{thm1}. This means that the difference between these two probabilities is negligible when compared with each of them, or what is the same thing, the conditional probability
\beqn\label{Cnd_Prb}
P_x[\sigma_A>n\, |\, \sigma_{\{0\}}>n, S_n =y]
\eeqn
approaches unity for every finite set $A$.

 The result in case  $xy<0$ and $|x|\wedge |y| \to\infty$ may be in a sense more interesting  and worthy  of note.
 Proposition \ref{prop4.1} in Section 6 entails  that under $x\vee |y| \leq M\sqrt n$,
 $$p_{\{0\}}^n(x,y) -p_A^n(x,y) = O((|x|\vee |y|)n^{-3/2}) \quad\mbox{ if}\;\; xy<0$$
   Combined with Theorem \ref{thm2} and  (a2) above this shows that  as $ x\wedge (-y) \to  \infty$ 
\beqn\label{DIFF}
\frac{p_A^n(x,y)}{p_{\{0\}}^n(x,y)}\; \longrightarrow\;  \left\{\begin{array} {ll} C^+_A/ C^+ \quad & \mbox{if}\quad E[|X|^3; X<0] <\infty,   \\
1 \quad & \mbox{if}\quad E[|X|^3; X<0] =\infty,
\end{array} \right.
\eeqn
thus the asymptotic form  of $p_A^n(x,y)$ in this regime does or does not  depend on $A$  according as $|x|^3p(x)$ is summable on $x<0$ or not, a  consequence suggested by the result in \cite{U1dm} mentioned just before Theorem \ref{thm2} but not fully expected to be true. One may wonder how
 the difference exhibited in (\ref{DIFF})   arises. This time the origin lies between $x$ and $y$ so that the walk must jump over the origin to avoid it whether the jump is short or long. 
 Now the difference in question is understood to accord with   different  behaviour of the walk:  in the first case  the walk that  reaches $y$  gets  close to the origin,  yet avoids it, whereas  in the second case it  clears the origin by a very long jump so as to  avoid $A$ simultaneously (see Corollary \ref{cor8} for more definite formulation).

\v2
We conclude this section by describing the main steps  of derivation of the  formula of Theorem \ref{thm1} (i)
 restricted to the case $|x|\vee |y| =o(\sqrt n)$.  
As mentioned before, the proof rests on  formula (\ref{crucial})
and its counter-relation
$$ P_x[\tau_{U(R)}  = \sigma_{(-\infty, -R]}<\sigma_A] =R^{-1} g^-_A(x)(1+o(1).$$
  In \cite{U1dm} we have shown that  (\ref{K}) holds uniformly in the regime $|x|\vee |y| = o(\sqrt n\,)$, which in particular  entails that as $|x|\wedge |y| \wedge n\to\infty$ under the condition  $|x|\vee |y| =o(\sqrt n\,)$
$$p_{\{0\}}^n(x,y) = \frac{2xy}{\sigma^2 n}p^n(y- x) + o\bigg(\frac{ |x| \vee |y|}{n^{3/2} }\bigg) \quad \mbox{if} \quad xy >0$$
 and
$$p_{\{0\}}^n(x,y) = o\bigg(\frac{ |x| \vee |y|}{n^{3/2} }\bigg)\quad \mbox{if} \quad xy < 0;$$
moreover, we shall see $p_{A}^n(x,y) \sim p_{\{0\}}^n(x,y)$ in the same limit scheme. 
 Taking $R=R_n=o(\sqrt n\,)$, so that $\tau_{U(R)} = o(n)$ (a.s.), we apply the strong Markov property at $\tau_{U(R)}$ and put all the relations mentioned above together to deduce  that for $|x| <|y| =o(\sqrt n)$, as $|y|\wedge n\to\infty$
$$p_A^n(x,y) = \frac{g_A^{{\rm sign}(y)}(x) |y|}{n}p^n(y- x) + o\bigg(\frac{(|x|+1) \vee |y|}{n^{3/2} }\bigg),$$
where $ {\rm sign}(y) $ is $+$ or $-$ according as $y$ is positive or negative. 
Finally  on applying the last relation to the time-reversed walk  the same  argument leading to it shows the 
formula of the theorem for $x$ and $y$ fixed.

In Section 2  we provide certain potential theoretic  facts  that  are used throughout the paper.  In Section 3.1 various  relations between $g_A$,  $u_A$ and $g^\pm_A$ are obtained, especially (\ref{1*})     is proved. In Section 3.2 we evaluate the probability of the walk exiting $U(R)$ without hitting $A$ for large $R$ and thereby prove (\ref{crucial}).  In Section 4 we state a few  known facts on   $p_{\{0\}}^n(x,y)$  and   $p^n_{(-\infty 0]}(x,y)$ and prove  some  results concerning them.   Proof of Theorem \ref{thm1} is given in Section 5.  In Section 6 the estimation of $p^n_A(x,y)$ is carried out   in case $xy<0$ and $|x|\wedge |y| \to\infty$, and in particular  the upper and lower bounds in Remark 2 (a) as well as 
Theorem  \ref{thm2}    are proved.  We provide in Section 7 some auxiliary results that are used in Remark 2 and  Sections  3, 5 and  6.

\section{Preliminaries  from the theory of 
one-dimensional random walks} 

In this section we review some potential theoretic results concerning one-dimensional random walks on $\Z$ that
are relevant to the present problem. The most of them are taken from Spitzer's book \cite{S}. Some of the  results and the
arguments taken from \cite{U1dm}  are included to make the description streamlined.   We shall designate by $C, C', C_1, \ldots$ etc. constants depending only on $p$ whose exact values are not significant for the present purpose and  may vary  at different occurrences of them. The letters  $x, y$ and $z$ are used to denote the  integers representing states of the walk. 

\v2

{\bf Potential function.}
 It is shown in \cite[Theorem 29.2]{S}  that 
$a(x+1) -a(x) \to  \pm 1/\sigma^2$ as $ x\to \pm\infty$,
which  implies   
$$a(x+z) -a(x) = \pm \frac{z}{\sigma^2}\{1+o(1)\} $$
with  $o(1) \to 0$ as $x\to \pm \infty$ uniformly for  $z$ with  $(x+z)x>0$,  and
 $a(x)/| x| \to 1/\sigma^2$
as $|x|\to  \infty$. It also holds that for all $x, y\in \Z$,
$$\sigma^2 a(x)\geq |x|$$ 
and 
$$\sum_{z\in \Z} p(z-x)a(z-y) = a^\dagger(x-y),$$
 that $a(x)$ is linear for $x\geq 0$ (resp. $x\leq 0$) if and only if the walk is left-continuous (resp. right-continuous, i.e., $P[X\geq 2] =0$) and that 
\beqn\label{a_rl}
\left\{ \begin{array}{ll} \sigma^2a(x)>x\; \;\; \mbox{for all} \; x>0 \quad & \mbox{ if  not left-continuous}, \\
\sigma^2a(-x)>- x\; \mbox{ for all}\;\;  x<0 \quad &\mbox{ if  not right-continuous.}
\end{array}\right.
\eeqn 
(Cf.  Theorem 28.1, Proposition 31.1 and  Proposition 30.3 of \cite{S} except for the strict inequality (\ref{a_rl})   which is found in \cite{U1dm}:see also (\ref{strct}), (\ref{A.03}) below). In what follows these relations will be  used frequently and not  noticed of their use. 

The results given in the rest of this section 
 will be used   in the proofs of Theorem \ref{thm4}  given in  Section 6 but not needed for  Theorem \ref{thm1} at least explicitly---some may be used in the proof of Theorem A  cited in Section 4.

\v2

{\bf Green's function on $[1,\infty)$.}
 Here we consider the walk  killed when it enters $(-\infty,0]$.  
 For $x=1,2,\ldots$,  let $v_{+}(x)$ (resp. $v_{-}(x)$) be   the probability that the strictly ascending  (resp.
descending) ladder process starting at $+1$ (resp. $-1)$ visits $x$ (resp $-x$): for $x\geq 2$
$$v_{-}(x)=P_{-1}[S_{\sigma_{(-\infty,-x]}}=-x]\quad and 
\quad v_{+}(x) = P_1[S_{\sigma_{[x,\infty)}}=x],
$$
and $v_{+}(1)=v_{-}(1)=1$, and  define 
 $$f^+(x) = E_0[\,|S_{\sigma_{(-\infty,-1]}}|\,](v_{-}(1)+\cdots +v_{-}(x)) $$
 and
 $$ f^-(x) = E_0[\,S_{\sigma_{[1, +\infty)}}\,](v_{+}(1)+\cdots +v_{+}(x)). $$
Then by the renewal theorem  $\lim_{x\to\infty} v_{-}(x) =1/E_0[\,|S_{\sigma_{(-\infty,-1]}}|\,]$ and similarly for $v_+(x)$, so that
 $f^+(x)\sim f^-(x) \sim x$ as $x\to\infty$. Under this boundary condition at $+\infty$,  $f^+(x)$ (resp. $f^+(x)$)  is the unique harmonic function with respect to  the walk  $S_n$ (resp $-S_n$) killed on $(-\infty, 0]$: 
\beqn\label{f_def}
~~f^{\pm}(x)=E[f^{\pm}(x\pm X);\, x\pm X >0]~~(x\geq 1)~~~\mbox{and}~~~ \lim_{x\to\infty} f^{\pm}(x)/x=1,
\eeqn
  (cf. Sections 18 and 19 of \cite{S}; in particular  Proposition 19.5 for the uniqueness).    (It is warned  that it is not  $[1,\infty)$  but   $[0,\infty)$  on which  the harmonic function is considered in \cite{S}.) It also holds \cite[Proposition 7.3]{U1dm}  that
 \beqn\label{strct}
 x< f^+(x) < \sigma^2a(x) \qquad \mbox{for}\quad x\geq 1
 \eeqn
unless   the walk is left continuous  (when $f^+(x) = \sigma^2 a(x) = x$ for $x>0$).
The Green function  $g_{\,(-\infty,0]}$ is expressed as
\beqn\label{g2}
g_{\,(-\infty,0]}(x,y)= \frac{2}{\sigma^2}\sum_{z=0}^{x\wedge y -1} v_+(x-z)v_-(y-z)~~~~~ (x,y >0).
\eeqn

 \v2
{\bf Hitting distribution of $(-\infty, 0]$.}
 Let   $H^x_{(-\infty,0]}(y)$ ($y\leq 0$) denote the  hitting distribution of $(-\infty,0]$ for the walk $S_n$ started at $x$. By the last exit decomposition  
 \beqn\label{h^-}
H^x_{(-\infty,0]}(y)=\sum_{w=1}^\infty g_{(-\infty,0]}(x,w)p(y-w).
\eeqn
As mentioned in Introduction   $H_{(-\infty, 0]}^{+\infty}(y) = \lim_{x\to\infty} H^x_{(-\infty,0]}(y)$
  is a probability.  By (\ref{g2})
\beqn\label{q}
H_{(-\infty, 0]}^{+\infty}(y)=\frac2{\sigma^2}E[f^-(y-X);X<y]=\frac{2}{\sigma^{2}}\sum_{w=1}^\infty f^-(w)p(y-w)~~~~~ (y\leq 0).
\eeqn
By summation by parts  we have
\beqn\label{Haa}
 H_{(-\infty, 0]}^{+\infty}(y)= \frac{ 2}{\sigma^{2}}  \sum_{w=1}^\infty v_-(w)F(y-w),
 \eeqn
 where  $F$ denotes the distribution function of $X$:
$$F(t) =P[X\leq t].$$
From (\ref{Haa}) it follows that  $H_{(-\infty, 0]}^{+\infty}(y)$ is monotone and   
$$H_{(-\infty,0]}^{+\infty}(y)\asymp \sum_{z<y}F(z); $$
 in particular 
\beqn\label{HiffE}
\sum H_{(-\infty,0]}^{+\infty}(y)|y|<\infty\quad \mbox{if and only if}\quad E[|X|^3; X<0]<\infty.
\eeqn
 In view of  (\ref{g2}) $g_{(-\infty,0]}(x,y)\leq Cx\wedge y\;\;(x,y>0) $, hence by $f^-(w)\leq C w$  and (\ref{q})
  \beqn\label{17}
   H^x_{(-\infty,0]}(y)\leq CH^{+\infty}_{(-\infty,0]}(y) \quad (y\leq 0);
\eeqn
also
\beqn\label{zF}
 H^{x}_{(-\infty,0]}(y) \leq Cx\sum_{w>0}  p(y-w) \leq C x F(z) \qquad(y<0).
 \eeqn 
\v2

{\bf  Green's function  on $\Z\setminus \{0\}$.}  This  is given by
\begin{eqnarray}
\label{greenft}
g(x,y)&:=&g_{\{0\}}(x,y) -\de_{x,0}  \nonumber\\
&=&a(x)+a(-y) -a(x-y)
\end{eqnarray}
(\cite[Proposition 29.4]{S}). Since $g(x,y) \leq g(x,x)\wedge g(y,y)$,
\beqn\label{zF2}
g(x,y) \leq C|x|\wedge |y|\quad \mbox{for all  $x$  and $y$} .
\eeqn
By the strong Markov property of the walk $S$
 we have 
\beqn\label{A.00}
\sum_{z=-\infty}^{0} H^{x}_{(-\infty,0]}(z) g(z,y) = g(x,y)\quad \mbox{for}\quad  y\leq 0 <x.
\eeqn
  By noting $g(\cdot, y)$ is bounded for  each fixed  $y$ we let  $x\to\infty$  to obtain
\beqn\label{A.000}
\sum_{z=-\infty}^{0} H^{+\infty}_{(-\infty,0]}(z) g(z,y) = a(-y) +\frac{y}{\sigma^{2}},
\eeqn
while on   letting $y\to -\infty$ in (\ref{A.00}) with the help of  (\ref{zF}) and (\ref{zF2})
\beqn\label{A.03}
\sum_{z=-\infty}^{0} H^{x}_{(-\infty,0]}(z) \bigg(a(z)- \frac{z}{\sigma^2}\bigg) = a(x)-\frac{x}{\sigma^2}\quad \mbox{for}\quad  x> 0.
\eeqn
In view of (\ref{17}) and (\ref{HiffE}) this shows (\ref{C_+}) as well as  the result stated following it, in particular 
 $\lim_{x\to\infty}[a(x)-x/\sigma^2]$ is finite if $E[|X|^3; X<0] <\infty$. 
 
\begin{lem}\label{lem3.1}\, Suppose $E[|X|^3; X<0]=\infty$.  Then,  as $x\to\infty$
\begin{eqnarray}\label{Hbb}
 a(x)- \frac{x}{\sigma^2} \, &\sim&\, \frac{ 2}{\sigma^{2}} \sum_{z=-x}^{-1} \sum_{w=-\infty}^{z} H_{(-\infty, 0]}^{+\infty}(w)\\
 &\sim&    \frac{4}{\sigma^4} \sum^{-1}_{z=-x-1} \sum^{z}_{w=-\infty} \sum^{w}_{j=-\infty} F(j).
 \label{Lamb0}
\end{eqnarray}
\end{lem}
\v2\n
\pf\,  
Let  $y:=-x \to  -\infty$ in (\ref{A.000}). Noting 
$$g(z, -x) = \left\{ \begin{array} {ll} 
[a(z)- \sigma^{-2}z] \{1+o(1)\} \sim 2\sigma^{-2}|z|  \quad & (-x<z<0),  \\[1mm]
 g(-x,-x)\{1+o(1)\} \sim 2\sigma^{-2}x \quad & (z\leq -x),
 \end{array}\right.
 $$
where $o(1)$ is uniform in  $z$, we then deduce that if $E[|X|^3; X<0]=\infty$,
\beqn\label{Hab}
a(x)- \frac{x}{\sigma^2}\, \sim\,  \frac2{\sigma^2}\sum_{z=-\infty}^{-1} (|z|\wedge x)H_{(-\infty, 0]}^{+\infty}(z). 
\eeqn
By summation by parts the right side equals that of (\ref{Hbb}) and 
(\ref {Lamb0}) follows  immediately from (\ref{Haa}) and (\ref{Hbb}).  The proof is  complete. 
\qed

\section{ The Green function and escape from $A$}
 
\subsection{Functions $u_A$, $g_A^+$ and $g_A^-$ }
 In analogy of  the corresponding formula for two-dimensional Brownian  motion (cf. \cite{H}) we define  $u_A(x)$, $x\in \Z$ by
\beqn\label{2_1}
u_A(x) = g_A(x,y) + a(x-y) - E_x[ a(S_{\sigma(A)}-y)].
\eeqn
 This conforms to   the definition in Section 1 in view of the next result proved in  \cite[Lemma 2.8]{Uf.s}.
\begin{lem}\label{lem2.1}\,The right side of (\ref{2_1}) is independent of $y\in \Z^2$ for all  $x\in \Z^2$. In particular for any $\xi_0\in A$
\beqn\label{2_2}
u_A(x)=  a^\dagger(x-\xi_0) - E_{x}[ a(S_{\sigma(A)}-\xi_0)].
\eeqn
\end{lem}
\v2
{\sc Remark 3.}  In  \cite{Uf.s} it is shown for a two-dimensional walk that  the right side  of (\ref{2_2})  is   dual-harmonic (i.e. invariant under the transform $h(y)  \mapsto \sum_z p(y-z)h(z)$) as a function of $y\in  \Z$,  hence   
does not depend on $y$   due to the fact that   
  the only bounded harmonic functions are constant functions \cite[Theorem 24.1]{S}, \cite[Proposition 5-20]{KSK}. The proof   applies without any change to every irreducible  recurrent random  walk (aperiodic or not) of any dimension. A non-lattice analogue of (\ref{2_1}) is obtained in \cite{PS}.

  For $x$ restricted on $A$,  (\ref{2_1})  reduces to the dual of the formula of Proposition 30.1 in  \cite{S}, where the dual of $u_A(x)$ (denoted by $\mu_A(x)$ therein)  is defined as the limit of $\sum_{z\in \Z^2} p^n(z-y)P_z[S_{\sigma(A)} =x]$ as $n\to \infty$. For $x\notin A$,  on the other hand,   an equivalent to (\ref{2_1}) (in a sense)   is found  in \cite[Proposition 4.6.3]{LL} where  $A$ may be  an infinite   set while  $p$ is  assumed to be of  finite range and symmetric. It is warned  that our $g_A(x,y)$ is different from   Green's function defined  in \cite{S} and \cite{LL} (see also (\ref{G-g})).

\v2\v2
Passing to the limit in (\ref{2_1}) we obtain
\beq
u_A(x) &=& \lim_{y\to +\infty} g_A(x,y) - (x- E_x[S_{\sigma(A)}])/\sigma^2\\
&=& g_A^+(x) - w_A(x),
\eeq
and  similarly $u_A(x) = g^-_A(x) + w_A(x)$, showing (\ref{1*}).    Substitution from  (\ref{2_2}) then yields  
expressions of $g_A^+$ and  $g_A^-$, which  we write down as
\beqn\label{repr_g}
g^+_A(x) = a^\dagger(x-\xi_0) +\frac{x}{\sigma^2} - E_x\Big[a(S_{\sigma(A)} -\xi_0) +\frac1{\sigma^2}S_{\sigma(A)}\Big]  
\eeqn
and
\beqn\label{repr_g-}
g^-_A(x) =  a^\dagger(x-\xi_0) -\frac{x}{\sigma^2} - E_x\Big[a(S_{\sigma(A)} -\xi_0) - \frac1{\sigma^2}S_{\sigma(A)}\Big], 
\eeqn
where $\xi_0$ is a point arbitrarily chosen  from $A$.

As $y\to\infty$ we have
$ a(S_{\sigma(A)}-y) = a(-y) - S_{\sigma(A)} +o(1),$
whence 
\begin{eqnarray}\label{2_3}
g_A(x,y)- g^+_A(x) &=& - w_A(x) - a(x-y) + E_x[a(S_{\sigma(A)}-y)]  \nonumber \\
&=& a(-y) - a(x-y) - \frac{x}{\sigma^2}  +o(1), 
\end{eqnarray}
where $o(1)\to 0$ as $y\to +\infty$ uniformly in $x \in \Z$. 
In Section 7  ((\ref{A.4}))  we shall see  that for $|x|\leq y$,
\beqn\label{A_2}
 \Big| a(-y) - a(x-y) - \frac{x}{\sigma^2} \Big | \leq    C\Big(a(-|x|) - \frac{|x|}{\sigma^2}\Big);
 \eeqn 
  in particular  the right side of (\ref{2_3})  is $(1+|x|)\times o(1)$ as  $y\to\infty$ uniformly for $|x|\leq y$.
This bound is presented here because of its bearing close relevance to (\ref{2_3}), although we shall use it only to show Lemma \ref{lem2.2} given at the end of this subsection and before proceeding to it we introduce some results   that are used in the proof of  (\ref{A.4})  as well as Lemma \ref{lem2.2}. 

Put
$$r^+_A = 1+ \max A\quad\mbox{and}\quad r^-_A =-1+  \min A,$$
so that  $A \subset [r_A^- +1, r_A^+ -1]$, and bring in  the following subsets of $\Z$:
$$V^+= \{x:  \exists n\geq 1, p^n_A(x, r^+_A )>0\},
\quad V^-= \{x\in \Z: \exists n\geq 1, p^n_A(x, r^-_A)>0\};$$
$$\widehat V^+= \{y : \exists n\geq 1, p^n_A(r^+_A, y)>0\},\quad \widehat V^-= \{y\in \Z: \exists n\geq 1, p^n_A(r^-_A,y)>0\};$$
in other words,  $V^+$ is the set of those points $x$ from which the walk  can enter  into  $[r_A^+, \infty)$  with a positive
probability, and similarly for $V^-_A$ and  $\widehat V_A^\pm$.  

Clearly  $g_A^+(x)>0$ for $x\geq r_A^+$.
By (\ref{harm_g}) 
\beqn\label{harm_g0}
g^+_A(x) = \sum_{z\in \Z} p^n_A(x,z) g^+_A(z)   \quad\mbox{for all}\quad  x\in \Z,\; n\geq 1;
\eeqn
hence  
$g^+_A(x) >0 $ for $x\in V^+$, while if   $x\notin V^+$, then $g_A(x,y) =0$ for $y\geq r_A^+$, hence $g^+(x)=0$; and similarly  for  $g^-_A$. Thus 
 \beqn\label{V}
 V^\pm = \{x: g_A^\pm(x)>0\}\quad \mbox{ and} \quad  \widehat V^\pm = \{y: g_{-A}^{\mp}(-y)>0\}.
 \eeqn

Write  $H^x_B(z)$ for $P_x[S_{\sigma(B)}=z]$ for $B\subset \Z$. By the strong Markov property
$$g_A(x,y) = \sum_{z\notin A,\, z\leq  r^+_A} H^{x}_{(-\infty,\,r^+_A]}(z)g_A(z,y),\quad   \; x>r^+_A,\; y< r^-_A,$$  
and, on  letting $y\to -\infty$,
$$g^-_A(x) = \sum_{z\notin A,\, z\leq  r^+_A} H^{x}_{(-\infty,\,r^+_A]}(z)g^-_A(z)$$  
as in the same way for  (\ref{A.03}).
Further let $x\to \infty$ to see that
$$C_A^+/\sigma^2 =\lim_{x\to \infty } g_{A}^-(x)=  \sum_{z\notin A,\, z\leq  r^+_A} H^{+\infty}_{(-\infty,\,r^+_A]}(z)g^-_A(z) \,(\leq \infty),$$
where the first equality is by definition and the second follows by (\ref{17}) i.e.,  $H^x_{(-\infty,0]}(z) \leq CH^\infty_{(-\infty,0]}(z)$ when  the last sum  is finite and by  Fatou's lemma when it  is infinity.
Taking the successive  limits in the reverse order  we can see that $\hat g_{A}^+(y) = g_{-A}^-(-y)$ converges to the same sum as $y\to -\infty$ (the justification is slightly  different but the result is obtained in an another way in Lemma \ref{lem7.1}  and we do not give any detail of it). As a consequence we find 
\beqn\label{2.26}
C_A^+ = C_{-A}^+= \sigma^2\sum_{z\notin A,\, z\leq r^+_A}  H^{+\infty}_{(-\infty,\,r^+_A]}(z)g^-_A(z).
\eeqn
Plainly   $H^{+\infty}_{(-\infty,\,r^+_A]}(r_A^+)>0$ and  these identities together show  

\begin{lem}\label{lem2.3} $C_A^+ =C_{-A}^+$ and the following three conditions are equivalent one another:
$$(a) \quad C_A^+>0; \quad (b) \quad g^-_A(r_A^+)>0; \quad (c)\quad     g^-_{-A}(-r_A^-) >0.$$
\end{lem}
\v2\n
\begin{lem}\label{lem2.4}  There exists a positive constant $c$ such that
\beqn\label{gg}
g_A^+(x)g^-_{-A}(-y) + g_A^-(x)g^+_{-A}(-y)  > c(1+|x|+|y|)
\eeqn
for all  $x, y$ for which  the left side is positive.
\end{lem}
\v2\n
\pf\; If  both $x$ and $y$ remain in a finite set, the assertion is trivial.   For reason of symmetry 
 we suppose $x>r_A^+\vee 1$.  Further we may suppose $y<-r_A^-$,  otherwise the left side of  (\ref{gg}) being asymptotic to $4xy/\sigma^4$.  If $g_A^-(r_A^+)=0$, then by Lemma \ref{lem2.3}  $g_{-A}^-(-r_A^-)=0$, hence the left side of (\ref{gg}) vanishes whenever $x>r_A^+$ and $y<-r_A^-$. It remains to  consider the case $g_A^-(r_A^+)>0$. In this case  by Lemma \ref{lem2.3}   $g_{-A}^-(-y)\geq  c_1$ ($y<-r_A^-$) for some $c_1$ so that  $g_A^+(x)g_{-A}^-(-y) \geq c_1'x$.  Similarly   $g_A^-(x)g_{-A}^+(-y) \geq c_2'|y|$. \qed
 
 \v2

\begin{lem}\label{lem2.2}  For each $M>1$,  
  uniformly for $-M<x\leq y$,  as $y\to +\infty$
$$g_A(x,y) = g_A^+(x)\{1+o(1)\};$$
and   for $x>0$ and  $y > \max A$, 
 \beqn\label{corA30}
 g_A(-x,y) \leq g_A^+(-x)\{1+o_y(1)\},
 \eeqn
where $o_y(1)$ is bounded and, as $y\to\infty$,  tends to zero   uniformly in $x$; in  (\ref{corA30})
the equality holds if and only if $E[ X^3; X>0]<\infty$.
\end{lem}
\v2\n
\pf\;  $g_A(x,y)$ is positive  for all $y\geq r_A^+$ or zero  for all $y\geq r_A^+$ according as $g_A^+(x) >0$ or $0$. Since  $g^+_A(x)\sim 2x/\sigma^2$ as $x\to\infty$,    (\ref{2_3}) and (\ref{A_2}) show the  first relation. 

 As for the second one  we may suppose that $0\in A$ and $g^+_A(r_A^-)>0$, and  let $\xi_0=0$ in (\ref{repr_g}). According to (the dual assertion of) Lemma \ref{lem2.3}, it follows from the latter condition that  $g_A^+(-x)$ is bounded away from zero for $x\geq  r_A^-$. Put  $\widehat \la(x) =  a(-x) -x/\sigma^2$. Then  (\ref{A_2}) becomes
 $$
 g_A(-x,y) - g_A^+(-x) = - \,\widehat\la(x+y) +\widehat\la(y) +o(1),
 $$
where $o(1)\to 0$ as  $y\to\infty$ uniformly in $x\in \Z$, and   since   $g^+_A(-x) = \widehat \la(x)+ O(1)$,  an application of (\ref{A.2}) concludes the second relation.  
 If $E[X^3; X>0] <\infty$ then  $\sigma^2 a(y)-y$ converges to a finite limit as $y\to\infty$, and by (\ref{2_3})
 it is easy to see that 
 $g_A(-x,y) - g_A^+(-x)$ tends to zero  as $x\wedge y\to\infty$ so that in (\ref{corA30}) the inequality sign may be replaced by the equality sign, whereas if 
$E[X^3; X>0] =\infty$, $g_A(-x,y)/g^+_A(-x)$ is made arbitrarily small by choosing an $x$ large enough.
The proof is complete.    \qed
\v2

\subsection{Probabilities of escape from $A$ and  an overshoot estimate }

What are advanced in this subsection are modifications of the corresponding results for the case $A=\{0\}$ 
that are given in \cite[Section 2.3]{U1dm} and the arguments are parallel  to those in it.

We have the identity
$$P_x[\sigma_{\{y\}}<\sigma_A] = \frac{g_A(x,y)}{g_A(y,y)}$$
whenever $y \neq x$. From (\ref{2_1}) and  (\ref{2_2}) we deduce 
$$g_A(y,y) = a(y)+a(-y) + O(1).$$
 Hence, 
from Lemma \ref{lem2.2}  it follows  that  uniformly for $-M<x< y$, 
\beqn\label{2_4}
P_x[\sigma_{\{y\}}<\sigma_A]  = \frac{\sigma^2g_A^+(x)}{2y}\{1+o(1)\} \quad  \mbox{as}\quad y\to +\infty
\eeqn
and similarly for the case when  $y\to -\infty$ and  $y<x < M$ with $g^-_A(x)$ replacing  $g^+_A(x)$ on the right side. In what follows the letter $R$ will always denote a positive integer. 

\begin{Prop}\label{prop2.1} \,  Uniformly for $ x< R$, 
$$P_x[\sigma_{[R,\infty)} < \sigma_A] =P_x[\sigma_{\{R\}}<\sigma_A] \{1+o(1)\}  \quad  \mbox{as}\quad R\to +\infty.$$
And furthermore,  for each  $M>1$,  as $R\to\infty$,  uniformly for $-M< x< R$,
 \beqn\label{2_41}
 P_x[\sigma_{[R,\infty)} < \sigma_A]  =\frac{\sigma^2g_A^+(x)}{2R}\{1+o(1)\}
 \eeqn
 and 
 \beqn\label{2_42}
 P_{x}[\sigma_{[R,\infty)} < \sigma_A] \leq \frac{\sigma^2 g_A^+(x)}{2R}\{1+o(1)\} \quad \mbox{uniformly for $ x< 0$}.
 \eeqn
\end{Prop}
\v2\n
\pf\, The difference $P_x[\sigma_{[R,\infty)}  <\sigma_A] - P_x[\sigma_{\{R\}}<\sigma_A]$ is expressed as 
\beqn\label{eq-1}
\sum_{z>R} P_x[\sigma_{\,[R,\infty)}< \sigma_A,\, S_{\sigma([R,\infty))} =z]P_z[ \sigma_{\{R\}}> \sigma_A].
\eeqn
For simplicity  suppose  $A \subset (-\infty,0]$. Then  this sum is   dominated by 
$$P_x[\sigma_{\,[R,\infty)}< \sigma_A] \sup_{z>R} P_z[\sigma_{\{R\}}> \sigma_{(-\infty,0]}],$$
of which   the supremum  tends to zero as $R\to\infty$  as is proved in   \cite[Lemma  2.3]{U1dm}.  Thus  we obtain the first relation of the proposition. 

The second relation  (\ref{2_41}) follows immediately from the first and  (\ref{2_4}). 
Similarly  (\ref{2_42}) follows by using   (\ref{corA30}), the second half of Lemma  \ref{lem2.2}
\qed

\begin{Prop}\label{prop2.2} (Overshoot estimate) ~~ For each $M>1$, uniformly for $-M\leq x<R$, as $R\to \infty$  
\beqn\label{ov_sh1}
\frac1{R}E_x[\,S_{\sigma([R,\infty))}-R\, |\, \sigma_{\,[R,\infty)}<\sigma_A\,]=o(1)
\eeqn
and 
 for  all $x\in \Z$ and  $\xi_0\in A$, 
\beqn\label{ov_sh0}E_x[\,S_{\sigma([R,\infty))};\, \sigma_{\,[R,\infty)}<\sigma_A\,]\leq  \frac12[\sigma^2 a(x-\xi_0) +x-\xi_0] +  c_A
\eeqn
for some constant  $c_A $  $(\leq \frac12 \sigma^2+ \frac12\sup_{\xi\in A}[\sigma^2 a(\xi -\xi_0) +\xi-\xi_0])$.   
\end{Prop}
\vskip2mm\n
{\it Proof.}\, Suppose $0\in A$ and $\xi_0=0$, which gives rise to no loss of generality since $R$ does not appear in the right side of (\ref{ov_sh0}).  (\ref{repr_g}) is then   reduced to 
$g_A^+(y) =  a^\dagger(y) + \sigma^{-2}y-  E_y[a(S_{\sigma(A)}) + \sigma^{-2}S_{\sigma(A)}] $.
That $g_A^+$ is non-negative and harmonic on $\Z\setminus A$ in the sense of the identity (\ref{harm_g})  implies that for  all $x\in \Z$,
$$E_x[g_A^+(S_{\sigma([R,\infty))}); \sigma_{\,[R,\infty)}<\sigma_A\,] \leq g^+_A(x),$$
which shows (\ref{ov_sh0}), because $y \leq \frac12(a(y) +\sigma^{-2}y)\leq \frac12[ \sigma^2 g^+_A(y) +  \sup_{\xi\in A}(\sigma^2 a(\xi) +\xi)]$ (for all $y\in \Z$) and
 $\sigma^2 g_A^+(x)\leq \sigma^2[1+a(x)] +x$.
For $R$ so large  that $g_A^+(z)$ 
 is monotone in $z\geq R$,   it also follows  that  
\begin{eqnarray} \label{ovsh}
0 &\leq& E_x[g_A^+(S_{\sigma([R,\infty))}) -g^+_A(R); \sigma_{\,[R,\infty)}<\sigma_A\,] \nonumber\\
 &\leq& g^+_A(x)-g^+_A(R)P_x[\sigma_{\,[R,\infty)}<\sigma_A\,].
\end{eqnarray}
Owing to   Proposition \ref{prop2.1}  the  last member  is dominated by
\[   g^+_A(x)\bigg(1- \frac{g^+_A(R)}{2R/\sigma^2}\bigg)+(1\vee x)\times o(1) = (1\vee x)\times o(1) ,
\]
where $o(1)\to 0$  as $R\to\infty$ uniformly for $-M\leq x<R$.   Dividing  by $RP_x[\sigma_{\,[R,\infty)}<\sigma_A]$ and applying Proposition \ref{prop2.1} again we find 
 $$\frac1{R}E_x[g_A^+(S_{\sigma([R,\infty))}) -g^+_A(R)| \sigma_{\,[R,\infty)}<\sigma_A\,] =  o(1).$$
Finally 
substitution from   
$R=\frac12 \sigma^2 g^+_A(R)\{1+o(1)\}$  ($R\to +\infty$)   yields the formula (\ref{ov_sh1}).
 \qed
 
 \v2
By Markov's inequality we deduce from Propositions  \ref{prop2.1} and  \ref{prop2.2} the following
\begin{Cor} \label{cor2.24}  Uniformly for $-M<x< R$ and for $R'>R$, as $R\to\infty$
\begin{eqnarray}\label{cor2.21}
 P_x[S_{\sigma([R,\infty))} \geq R',  \sigma_{[R,\infty)}<\sigma_A] &=& \frac{RP_x[ \sigma_{[R,\infty)}<\sigma_A]}{R'-R}\times o(1) \nonumber \\
 &=& \frac{g_A^+(x)}{R'-R} \times o(1).
 \end{eqnarray}
 \end{Cor} 
\begin{Cor} \label{cor20}  For any  $\e>0$, uniformly for $-M<x< R$, as $R\to\infty$ 
\beqn\label{eq_cor20}
 E_x[S_{\sigma([R,\infty))}; S_{\sigma([R,\infty))} \geq (1+\e)R,  \sigma_{[R,\infty)}<\sigma_A] 
  = g_A^+(x) \times o(1).
 \eeqn
 \end{Cor} 
 \v2\n
 \pf\, 
We have only to decompose  $S_{\sigma([R,\infty))}$ into the sum of  $S_{\sigma([R,\infty))}-R$ and  $R$ and to apply
Propositions \ref{prop2.2} and  \ref{prop2.1} for the first term and Corollary \ref{cor2.24} for the second.
\qed

\begin{Prop}\label{prop3}\,  For each $M >1$, as $R\to\infty$
\beq
 P_x[ \tau_{U(R)} = \sigma_{[R,\infty)} < \sigma_A] &=& P_x[\sigma_{[R,\infty)} < \sigma_A] \{1+o(1)\} \\
 &=& \frac{\sigma^2 g^+_A(x)}{2R}\{1+o(1)\},
 \eeq
where $o(1)\to 0$ uniformly for $-M\leq x\leq R$. 
\end{Prop}
\v2\n
{\it Proof.}\, 
Given  $\xi_0 \in A$, we put $f(z)= a(z-\xi_0)+ \sigma^{-2}z$, so that in view of (\ref{repr_g})
$$g^+_A(x) =\de_{\xi_0,x}+ f(x)- E_x[f(S_{\sigma_A})].$$ 
Suppose $x\neq \xi_0$,    the case $x=\xi_0$ being  readily reduced to this case (see (\ref{harm_g})). 
In view of  the  optional sampling theorem  $M_n=f(S_{n\wedge \sigma_A \wedge \tau_{U(R)}})$ is  then a martingale under the law $P_x$. It is uniformly integrable  and hence, according to  the martingale convergence theorem,  
$$f(x)=E_x[f(S_{\sigma_A \wedge \tau_{U(R)}} )].$$
 Breaking this expectation 
according as  $\tau_{U(R)}$ is larger  or smaller than $\sigma_A$ we find 
\beq
g_A^+(x)  &=& E_x[f(S_{\sigma_A \wedge \tau_{U(R)} })] - E_x[f(S_{\sigma_A})]\\
&=& E_x[ f(S_{\tau_{U(R)}})\,;\, \tau_{U(R)} < \sigma_A]
 - E_x[ f(S_{\sigma_{A}}); \tau_{U(R)} < \sigma_A].
\eeq
The last expectation   uniformly converges to zero as $R\to\infty$. In view of Proposition \ref{prop2.1} and the relation  $f(R) = 2R/\sigma^2\{1+o(1)\}$, for the proof of the proposition  it  therefore  suffices to show  
\beqn\label{p5}
 E_x[ f(S_{\tau_{U(R)}})\,;\, \tau_{U(R)} < \sigma_A]- f(R)P_x[ \sigma_{[R,\infty)} =\tau_{U(R)}  < \sigma_A]= (1\vee |x|) \times o(1)
 \eeqn
  as $R\to\infty$ uniformly for $-M<x<R$.
   By the dual relation of  (\ref{ov_sh0}) 
 we deduce 
 \beq
  E_x[ f(S_{\tau_{U(R)}})\,;\, \sigma_{(-\infty, -R]} = \tau_{U(R)} < \sigma_A]
 & \leq &
  \bigg[\sup_{y\leq -R}\frac{f(y)}{-y}\bigg]E_x[ |S_{\sigma_{(-\infty, -R]}}|\,;\, \sigma_{(-\infty, -R]} < \sigma_A]\\
&\leq&    [\sigma^2a (-x)-x +1] \times o(1)
  \eeq  
  (uniformly for all $x$).
On the other hand, owing to the overshoot estimate in  (\ref{ov_sh1}) 
 \beq
E_x[ f(S_{\tau_{U(R)}}) - f(R) \,;\, \sigma_{[R,\infty)} =\tau_{U(R)} < \sigma_A] = (x\vee 1) \times o(1)
\eeq
uniformly for $-M<x <R$. Adding these relations we obtain  (\ref{p5}).
\qed

\begin{Cor}\label{cor2.1}\,  For any  $M > 1$,  uniformly for  $-M<x\leq R$ and $z \geq R$, as $R\to\infty$
\begin{eqnarray}\label{A_6}
&& P_x[ S_{\sigma_{[R,\infty)}} =z \,|\,\tau_{U(R)}=  \sigma_{[R,\infty)}<\sigma_A] \nonumber\\
&&\quad  \leq 
P_x[ S_{\sigma_{[R,\infty)}} =z \,|\,  \sigma_{[R,\infty)}<\sigma_A]\{1+o(1)\}.
\end{eqnarray}
\end{Cor}
\v2\n
\pf\, Denote the conditional probabilities on the left and on  the right  of  (\ref{A_6}) by $\mu(z) $ and 
$\tilde \mu(z)$, respectively.
Then   we see 
\beq
\mu(z) &\leq& \frac{P_x[ \sigma_{[R,\infty)}<\sigma_A, S_{\sigma_{[R,\infty)}}=z]}{P_x[\tau_{U(R)}=  \sigma_{[R,\infty)}<\sigma_A] }  \\
&=& \tilde \mu(z)\frac{P_x[ \sigma_{[R,\infty)}<\sigma_A]}{P_x[\tau_{U(R)}=  \sigma_{[R,\infty)}<\sigma_A]}.
\eeq
According to   Proposition \ref{prop3} this yields  $\mu(z) \leq  \tilde \mu(z)\{1+o(1)\}$, as desired.  \qed
\v2

The next result concerns the probability of the walk   escaping  from  $A$, which, though not applied in this paper, we record here. 
\begin{Cor}\label{cor3}~~ Uniformly for $ |x|<R$, as $R\to\infty$
$$P_x[ \tau_{U(R)} <\sigma_A]= \frac{\sigma^2 u_A(x)}{R}\{1+o(1)\}.
$$
\end{Cor}
\vskip2mm\n
\pf~   The inclusion-exclusion formula  derives the assertion of the proposition  from Proposition \ref{prop3} if  
 $P_x[ \sigma_{(-\infty,-R]} \vee  \sigma_{[R,\infty)}<\sigma_{\{0\}}]=o({x}/{R})$
 uniformly for  $0<|x|<R$, which is readily verified (cf. \cite[Proposition 2.4]{U1dm}). ~~ \qed

\section{ Results for a single point set and  a half line}
Put
\beqn\label{sf_g}   \qquad \qquad {\sf g}_t(u) = \frac1{\sqrt{2\pi \sigma^2 t}}e^{- u^2/2\sigma^2 t}\qquad  (u\in \R, t>0).
\eeqn
We shall apply the following version of local limit theorem (cf. e.g., \cite{S}, \cite{Ws}):  as $n\to\infty$, uniformly for $x\in \Z$ with $p^n(x)>0$,
\beqn\label{llt}
 p^n(x) = \nu {\sf g}_n(x) + o\bigg(\frac{\sqrt n}{n\vee x^2}\bigg),
\eeqn
 where $\nu$ designates the (temporal) period of the walk: $\nu =$ g.c.d. of $\{n: p^n(0) >0\}$.

 \subsection{ Results for the case $A=\{0\}$}

Here we state two results from \cite{U1dm}  and \cite{Ufh} for the case $A=\{0\}$ that are used later. 

  \v2\n
{\bf Theorem A.} \, 
{\it  Given a constant $M>1$, the following asymptotic estimates of $p_{\{0\}}^n(x,y)$ as $n\to \infty$,  stated in  three cases of constraints on $ x$ and $y$, hold true  uniformly for $x$ and $y$ subject to the respective constraints. }

\v2
{\bf (i)}~ {\it Under  $|x|\vee|y|< M \sqrt n$ and $|x|\wedge |y|=o(\sqrt n)$,}
\beqn\label{(i)}
p_{\{0\}}^n(x,y)=\frac{\sigma^4a^\dagger(x)a(-y)+xy}{\sigma^2n}\,p^n(y-x)+o\bigg(\frac{(|x|\vee1)|y|}{n^{3/2}}\bigg).
\eeqn

\v2
{\bf (ii)}~  {\it Under  $M^{-1}\sqrt n < |x|,\,|y|< M \sqrt n$ (both $|x|$ and $|y|$ are between the two extremes),}
\begin{eqnarray}\label{q(ii)}
p_{\{0\}}^n(x,y)&=&\nu \Big[{\sf g}_{n}(y-x)-{\sf g}_{n}(y+x)\Big]+o\bigg(\frac{1}{\sqrt n}\bigg) \nonumber \\
&& \quad \quad\quad \quad \quad\quad \quad\quad \quad  \mbox{if} \quad xy>0\,\,\mbox{ and}\,\,  p^n(y-x) > 0,\\
p_{\{0\}}^n(x,y)&=&o\bigg(\frac{1}{\sqrt n}\bigg)  \,\,\quad \quad\quad \quad\quad\mbox{if}\quad\quad xy<0.
\label{(ii)}
\end{eqnarray}

\v2
{\bf (iii)}~ {\it  Let $1\leq  |x|\wedge |y| <\sqrt n< |x|\vee|y|$.  ~ Then, if $E|X|^{2+\de}<\infty$ for some $ \de\geq 0$, }
\beqn\label{eq4.6}
p_{\{0\}}^n(x,y)=O\bigg(\frac{|x|\wedge |y|}{|x|\vee|y|}{\sf g}_{4n}(|x|\vee|y|)\bigg)+\frac{|x|\wedge |y|}{(|x|\vee|y|)^{2+\de}}\times o(1).
 \eeqn

 \v2\n
{\bf Theorem B.} \,
{\it Uniformly in $x$, as $n\to\infty$}
$$
P_x[\sigma_{\{0\}} =n] =\frac{\sigma^2 a^\dagger(x)}{n}p^n(-x) + \frac{|x|\vee1}{n^{3/2} \vee |x|^3}\times o(1).
$$

Theorem A  is Theorem 1.1 
of \cite{U1dm}. Theorem B is an immediate consequence of Corollary 1.1 of \cite{Ufh} and the local limit theorem (\ref{llt}), if $\nu=1$, i.e., if the walk $S$ is  aperiodic. For $\nu>1$, apply this 
result to the walk
 $\tilde S_n = \nu^{-1} S_{\nu n}$, whose increment   has variance  $\sigma^2/\nu$.  
\v2
From Theorem A it follows  that as $|x|\wedge |y| \to\infty$ 
\beqn\label{eq4.7}
\left\{ \begin{array} {ll}
 {\rm (a)} \;\;  p_{\{0\}}^n(x,y)  \asymp |xy|n^{-3/2} \quad& (xy>0, |x|\vee |y| =O(\sqrt n\,)), \\[2mm]
 {\rm (b)} \;\; p_{\{0\}}^n(x,y)   = o(|xy|n^{-3/2}) \quad & (xy<0).
 \end{array} \right.
\eeqn
Indeed, (a) is obvious, whereas  for any  $M$, (b) follows  for $|x|\vee |y|< M\sqrt n$ by (i) on  one hand and  $p_{\{0\}}^n(x,y)< C |xy|n^{-3/2}M^{-2}$ for  $|x|\vee |y|\geq  M\sqrt n$ by  (ii) on the other hand. 
 \v2
 {\it Proof of Proposition \ref{prop1.1}.} 
The bound of the proposition  follows from the local limit theorem if $|x|\wedge |y|\geq \sqrt n$, from 
 (\ref{eq4.6})  if  $|x|\wedge |y| < \sqrt n \leq  |x|\vee |y|$,  and from   (\ref{(i)}) and (\ref{eq4.7})  if $|x|\vee |y| <\sqrt n$. \qed

\subsection{ Space-time  distribution of entrance into $(-\infty, 0]$}

 Here we consider the walk  killed when it enters $(-\infty,0]$.  The results given in this subsection 
 will be used only in  Section 6 where the case $xy<0$ is intensively studied.

The following result is  obtained as a special case  of    \cite[Proposition 11]{D} that concerns asymptotically stable  random walks (cf. also \cite[Theorem1.3]{U1dm}).

 \v2\n
{\bf Theorem C.}\, {\it For each  $M>1$, uniformly for $0< x, y\leq M\sqrt n$, as $xy/n\to 0$}
$$p_{(-\infty,0]}^n(x,y)=\frac{2f^+(x)f^-(y)}{\sigma^2n}p^n(y-x)\{1+o(1)\}.$$

In the regime   $M^-\sqrt n<x,y < M\sqrt n$ excluded from this theorem  $p_{(-\infty,0]}^n(x,y)$
has the same asymptotic form as $p_{\{0\}}^n(x)$ which is  given in Theorem A (ii), namely 
\beqn\label{q(ii)2}
p_{(-\infty,0]}^n(x,y)=\nu \Big[{\sf g}_{n}(y-x)-{\sf g}_{n}(y+x)\Big]+o\bigg(\frac{1}{\sqrt n}\bigg) 
\eeqn
provided $p^n(y-x) > 0$  (see Remark 4 after Lemma \ref{lem4.1} in the next section).

 \vskip2mm

 From Theorem C we  derive an asymptotic form of the space-time  distribution of  the first entrance into $(-\infty,0]$, which we denote by $h_x(n,y)$:  for $y\leq 0$
 $$h_x(n,y) =P_x[S_{\sigma_{(-\infty, 0]}}=y, \sigma_{(-\infty, 0]}=n].$$
We suppose $a$ as well as  $H^\infty_{(-\infty,0]}$  to be extended to continuous 
 variables by linear interpolation  for notational convenience.
 \v2

\begin{lem}\label{thm1.4}~  Define  $\a_n(x,y)$  for $y\leq 0<x $  via the equation
\beqn\label{h0}
h_x(n,y)=\frac{\nu f^+(x){\sf g}_n(x)}{n}\Big[ H_{(-\infty, 0]}^{+\infty}(y) +\a_n(x,y)\Big] \eeqn
if $p^n(y-x)>0$,  and put $\a_n(x,y)= 0$ if $p^n(y-x) =0$. 
Then for each $\e>0$ and $M\geq 1$, $\a^n(x,y)$ can be decomposed as
$$\a_n(x,y) = \b_{n,\e}(x,y) + H_{(-\infty, 0]}^{+\infty}(y) \times o_{\e,M}(1),$$
  where $o_{\e,M}(1)$ is bounded and,  as $n\to\infty$ and $\e\to 0$ in this order, approaches zero uniformly for $0<x< M\sqrt n$ and $y<0$, and the function of $y\leq 0$ defined by
 $$\b_{n,\e}(y) :=\sup_{n/2\leq k\leq n}\,  \sup_{0< x <M\sqrt n} |\b_{k,\e}(x,y)|,$$
satisfies   
\beqn\label{alpha1}
\b_{n,\e}(y) \leq c_MH_{(-\infty, 0]}^{+\infty}(y),  \quad  \lim_{n\to\infty} \sum_{y =-\infty}^0\b_{n,\e}(y)  = 0
\eeqn
 and 
 \beqn\label{alpha2}
      \sum_{z=-\infty }^{0}\b_{n,\e}(z)(|z|\vee |y|) \leq C|y|\frac{\sigma^2 a(\e \sqrt n) - \e \sqrt n}{\e \sqrt n}.\eeqn
\end{lem}
\v2\n
\v2\n
\pf\, 
Suppose the walk is aperiodic  for simplicity. We have  the representation
 \beqn\label{h00}
 h_x(n,y) = \sum_{w=1}^\infty  p_{(-\infty,0]}^{n-1}(x,w)p(y-w).
 \eeqn
In view of Theorem C and the local limit theorem (\ref{llt}), for each $\e>0$,  the sum on the right over $w< 4\e\sqrt n$ may be   replaced by
$$\frac{f^+(x) {\sf g}_n(x)}{n} \cdot  \frac{2}{\sigma^2}\sum_{1\leq w\leq 4\e\sqrt n}f^-(w)  p(y-w) (1+o_{\e,M}(1)).$$
By definition  $\a_n(x,y)$ is then expressed as 
\beqn\label{alpha4}
\a_n(x,y) = -\frac{2}{\sigma^2}\sum_{1\leq w\leq 4\e\sqrt n}f^-(w)  p(y-w) \times o_{\e,M}(1)
+\b_{n,\e}(x,y),
\eeqn
where $\b_{n,\e}(x, y)$ is defined to be
$$-\frac{2}{\sigma^2} \sum_{ w > 4\e\sqrt n}f^-(w)  p(y-w) + \bigg( \frac{f^+(x) {\sf g}_n(x)}{n}\bigg)^{-1}  \sum_{w>4\e \sqrt n}^\infty  p_{(-\infty,0]}^{n-1}(x,w)p(y-w).$$
By Proposition \ref{prop1.1}
 $p^{n-1}_{(-\infty,0]}(x,w) \leq x wn^{-3/2}$
and  we infer that
$$|\b_{n,\e}(x,y)|\leq  C\sum_{w>4\e\sqrt n}  wp(y-w),$$
where $C$ may depend on $M$.
 The first term on the right side of (\ref{alpha4}) may be replaced by  $H_{(-\infty, 0]}^{+\infty}(y) \times o_{\e,M}(1)$ and 
     the verification of  (\ref{alpha1}) is now immediate. 

For the proof of (\ref{alpha2})   observe  first that by the change  of  variable $w= 4\e\sqrt n+ z$
$$
|\b_{n,\e}(x,y)| \leq C_1 H_{(-\infty, 0]}^{+\infty}( y-4\e\sqrt n) + C_2\e\sqrt n F(y-4\e\sqrt n),$$
and then that on summing by parts and  abbreviating $H_{(-\infty, 0]}^{+\infty}$ to $H$,
\beq
\sum_{z=-\infty}^0 (|z|\wedge |y|) |\b_{n,\e}(x,z)| &=&\sum_{z=y}^{-1}  \sum_{w=-\infty}^z |\b_{n,\e}(x,w)|\\
&\leq& C_3\sum_{z=y}^{-1}\bigg[ \sum_{w=-\infty}^{z-4\e\sqrt n} H(w) + 4\e\sqrt n H(z-4\e\sqrt n)\bigg]\\ 
&\leq &C_3 |y|\bigg[ \sum_{w=-\infty}^{y-4\e\sqrt n} H(w) + 4\e\sqrt n H(-4\e\sqrt n)\bigg]
\eeq
where (\ref{Haa}) is applied for the first inequality and the monotonicity of $H(w)$ for the  second. Since $4\e\sqrt n H(-4\e\sqrt n) \leq 8 \sum_{w\leq -2\e\sqrt n}H(w)$ and since
$$\sum_{w\leq t} H(w) \leq C'' \sum_{w\leq t} H(w)g(w,t)/|t| \leq C''[a(-t)+t/\sigma^2]/|t|   \qquad (t<0), 
$$ we conclude  
$$\sum_{z=-\infty}^0 (|z|\wedge |y|) |\b_{n,\e}(x,z)|  \leq C|y| \frac{a(-2\e \sqrt n) -2\e \sqrt n \,}{2\e \sqrt n }.$$
This bound is valid with  $\b_{\e,n}(z)$ in place of  $|\b_{n,\e}(x,z)|$ as easily checked by following the    above derivation  with  $k$ instead of   $n$. Thus (\ref{alpha2}) is proved.
 \qed
\v2

From the cheaper estimate (\ref{alpha1}) we have the asymptotic form of the  distribution of hitting time of $(-\infty,0]$. We record it here for later citations: as $n\to\infty$
\beqn\label{hitt_HL}
\frac1{\nu}\sum_{j=0}^{\nu-1}P_x[\sigma_{(-\infty,0]} = n+j] \sim  \frac{ f^+(x){\sf g}_n(x)}{n}
\eeqn
 uniformly for $0<x< M\sqrt n$
  (cf. \cite[Corollary 1.1]{U1dm} for an estimate for $x>\sqrt n$).


\section{ Proof of Theorem \ref{thm1}}  
We break   Theorem \ref{thm1} into   three assertions by dividing   the range of variables into three regimes according as
 $|x|\wedge |y|\to\infty$,  $|x|\vee |y|\to\infty$ or  $|x|\vee |y| = O(1)$, of which  the first  case is dealt with by Lemma \ref{lem4.1}, the second by Lemma \ref{lem4.2} and the third  in the \lq {\it Proof of   Theorem \ref{thm1}}\,' given at  the end of this section. 
\v2\n
\v2

\begin {lem}\label{lem4.1}    As   $x\wedge y \wedge n\to\infty$ under the constraint
$ x\vee  y < M\sqrt n$,  
\beqn\label{3_1}
p_{\{0\}}^n(x,y) - p_A^n(x,y) = p_{\{0\}}^n(x,y)\times O\Big(\frac{1}{x\wedge y}\Big)
\eeqn
and \quad $p_A^n(-x,y) = o(xyn^{-3/2}).$
\end{lem}
\v2\n
\pf\, The second relation  is   (\ref{eq4.7}\,b). For the proof of (\ref{3_1}) 
suppose  $0\in A$ for simplicity---otherwise replace $\{0\}$ by $\{c\}$ with $c\in A$---so that 
\beqn\label{100} 
p_{\{0\}}^n( x,y) - p^n_A(x,y) = \sum_{k=0}^n \sum_{\xi\in A\setminus \{0\}} P_x[ \sigma_A=k, S_{k}=\xi] p_{\{0\}}^{n-k}(\xi, y).
\eeqn
We split the outer sum at $k=\lfloor n/2\rfloor$,  and denote the  double sum   restricted to $k<n/2$ by   $I_{[0,n/2)}$  and  the other part of the sum by $I_{[n/2,n]}$. According to Theorem A (i)
$$I_{[0, n/2)} \leq \sharp A \sup_{k<n/2}\sup_{\xi\in A}  p_{\{0\}}^{n-k}(\xi, y) \leq C \frac{y}{n}p^n(y),$$
where $\sharp A$ denotes the cardinality of $A$.
In a similar way, applying  Theorem B  to find  that 
$$I_{[n/2, n]} \leq \frac{Cx}{n^{3/2}}\sup_{\xi}  \sum_{n/2\leq k\leq n} p_{\{0\}}^{n-k}( \xi, y) \leq \frac{C'x}{n^{3/2}}.$$
where  it is used for the second inequality  that 
$\sum_{j=1}^\infty p_{\{0\}}^j(\xi,y)  = g_{\{0\}}(\xi,y) \leq C''$ ($\xi \leq \max A$, $y>0$). 
Adding  these two  bounds concludes the proof of Lemma \ref{lem4.1}. 
\qed
\v2

{\sc Remark 4.}\, With minor modification    the above proof also shows
\beqn\label{Rem4}
p^n_{\{0\}}(x,y) - p^n_{(-\infty,0]}(x,y) = p^n_{(-\infty,0]}(x,y)\times[O(1/y)+ o(\la(x)/x)].
\eeqn
Indeed 
the evaluation of $I_{[n/2,n]}$ made in it  applies to the  case $A=(-\infty,0]$ if we use  (\ref{hitt_HL}) instead of Theorem B, while 
for  $I_{[0,n/2]}$ we apply  (\ref{eq4.7}\,b), (\ref{A.03}) and Theorem C in turn to obtain
$$I_{[0,n/2]} 
=  \sum_{\xi\leq 0}H^x_{(-\infty,0]}(\xi)\times o\bigg(\frac{|\xi|y}{ n^{3/2}}\bigg)=o\bigg(\frac{\la(x)y}{n^{3/2}}\bigg) = p^n_{(-\infty,0]}(x,y)\times o\bigg(\frac{\la(x)}{x}\bigg). $$
Although (\ref{Rem4}) is not used in this paper (apart from (\ref{q(ii)2})), it indicates the possibility for another approach based on the results for the walks killed on  half lines.
\v2

 In view of Theorem A and the  local  limit theorem (\ref{llt})  it follows from Lemma \ref{lem4.1} that as   $x\wedge y \wedge n\to\infty$ under the constraint
$ x\vee  y < M\sqrt n$,  
\beqn\label{eq4.4}
p_A^n(x,y) =\nu [{\sf g}_n(y-x)- {\sf g}_n(x+y)](1+ o(1)),
\eeqn
provided  $p^n(y-x)>0$. This   proves the second half of Theorem \ref{thm1}.   
 \v2

\begin {lem}\label{lem4.2} \,  For any  $M\geq 1$,  uniformly for $-M <x < \sqrt n/\lg n$, as  $ y\to\infty$ under $y <M \sqrt n$,
$$p_A^n(x,y) = \frac{g_A^+(x)y}{n}p^n(y-x)\{1+o(1)\};$$
and   uniformly for $-M< y< \sqrt n/\lg n$, as $x \to +\infty$ under $x <M \sqrt n$,
$$p_A^n(x,y) = \frac{xg^+_{-A}(-y)}{n}p^n(y-x)\{1+o(1)\}.$$
\end{lem}

\v2\n
\pf\,  We have only to prove the first relation, the second one being  the dual of it. We may suppose $g_A^+(x)>0$, otherwise $p_A^n(x,y)=0$ for all large $y$. Put $R= \lfloor \sqrt n/ \lg n\rfloor$ and  $N = mR^2\lfloor \lg n\rfloor$ with  $m$ determined shortly and  decompose
\begin{eqnarray}\label{decomp}
p^n_A(x,y) &=& \sum_{k=1}^{N-1} \sum_{|z|\geq R}P_x[\tau_{U(R)}=k<\sigma_A, S_k=z]p_A^{n-k}(z,y) \\
&& +\,\, \e(x,y;R).\nonumber
\end{eqnarray}
Here
\beqn\label{e}
\e(x,y;R) = \sum_{z}P_x[ \tau_{U(R)\setminus A} \geq N, S_N =z]p_A^{n-N}(z,y)
\eeqn
and $\lfloor t\rfloor$ denotes the largest integer that does not exceed $t$.
Using  the central limit theorem we deduce  $\sup_{x\in U(R)}P_x[\tau_{U(R)} > R^2] < c$ for all sufficiently large $R$ with a universal constant $c<1$, and hence that  
$$\e(x,y;R)  \leq  Ce^{-\la N/R^2} /\sqrt n$$
with  $\la = -\lg c$. Now take $m=2/\la$ so that 
\beqn\label{N_R}
\e(x,y;R)  =  O(n^{-2}) ;
\eeqn
hence  $\e(x,y;R) $ is negligible. (Note that $p_A^n(x,y) =0$ if $g^+_A(x)=0$.)

As for the double sum in (\ref{decomp}) we first notice that   the contribution from the half line
$z\leq -R$ is negligible in comparison with that from $z\geq R$ owing to  (\ref{2_42}) and the second relation of Lemma \ref{lem4.1}.  
It  therefore remains to evaluate 
\beqn\label{decomp1}
W :=\sum_{k=1}^{N-1} \sum_{z\geq R} P_x[ \tau_{U(R)} = k <\sigma_A, S_k=z]p_{A}^{n-k}(z, y)
\eeqn
so as to verify
\beqn\label{d_sum}
W =\frac{g_A^+(x)y}{n}p^n(y-x)\{1+o(1)\}.
\eeqn
From  (\ref{3_1}) and Theorem A (i) and (ii)  of Section 4.1 it follows that if $k<N$,
\beqn\label{5.9_1}
p_{A}^{n-k}(z, y) \left\{ \begin{array}{ll} {\displaystyle  
= \frac{2zy}{\sigma^2n}p^n(y-z)\{1+ o(1)\} } \quad &\mbox{for}\quad R\leq  z< 2R, \\[3mm]
\leq  zy /n^{3/2} \quad &\mbox{for}\quad 2R \leq  z<\sqrt n. 
\end{array}\right.
 \eeqn
By the same reasoning as for (\ref{N_R}) it also follows that
\beq
 P_x[\sigma_{[R,\infty)} = \tau_{U(R)} < \sigma_A] =\sum_{k=1}^{N-1}  \sum_{z\geq R} P_x[ \tau_{U(R)} = k <\sigma_A, S_k=z] - O(n^{-2}).
\eeq
From these observations we infer that
\begin{eqnarray}\label{W-P}
&&\bigg| W - P_x[\sigma_{[R,\infty)} = \tau_{U(R)}  <\sigma_A]  \frac{2Ry}{\sigma^2n}p^n(y-x) \bigg| \nonumber\\
&&\leq \, \sum_{R\leq z< 2R} P_x[ \tau_{U(R)}  <\sigma_A, S_{\tau_{U(R)}}=z] \bigg|\frac{2zy}{\sigma^2n}p^n(y-z)\{1+o(1)\}  - \frac{2Ry}{\sigma^2 n}p^n(y-x)\bigg| \nonumber\\
 &&\quad \quad + \,  \frac{y}{n^{3/2}}E_x[ S_{\tau_{U(R)}} ; \tau_{U(R)} <\sigma_A,\, S_{\tau_{U(R)}} \geq 2R]    \nonumber \\
 &&\quad \quad +\,  \frac{CRy}{n^{3/2}} P_x[ \tau_{U(R)}  <\sigma_A, \,S_{\tau_{U(R)}}\geq 2R] \nonumber\\
&&\quad \quad + \, 
   C n^{-2}.
\end{eqnarray}
On the right side   we may replace  $\tau_{U(R)}$ by $\sigma_{[R,\infty)}$ for obvious reason.
In view of the  local  limit theorem,  $p^n(y-z) = p^n(y-x)\{1+o(1)\}$ 
uniformly for $|z|< 2R$,  hence 
$$|zp^n(y-z) - Rp^n(y-x)|\leq zp^n(y-x)\times o(1)+  |z-R|/\sqrt n,$$
 provided $p^n(y-x) p^n(y-z)>0$, and 
 we infer  that
the first  term  on the right side of \\
 (\ref{W-P})  is  at most a constant multiple of 
\beq
 \sum_{R\leq z< 2R} P_x[ \sigma_{[R,\infty)} <\sigma_A, S_k=z]\bigg[\frac{2zy}{\sigma^2n}p^n(y-x)\times o(1) +\frac{2 (z-R)y}{\sigma^2 n^{3/2}}\bigg],
\eeq
which, on applying   Propositions \ref{prop2.2} and \ref{prop2.1} in turn,  is 
$$ P_x[\sigma_{[R,\infty)}  <\sigma_A] \frac{Ry}{ n^{3/2}}\times o(1) =\frac{yg^+_A(x)}{n^{3/2}}\times o(1).$$
By Corollary \ref{cor20}  (in Section 3.2)  the second term on the right side of (\ref{W-P}) is at most
$$ \frac{yg^+_A(x)}{n^{3/2}}\times o(1).$$
Similarly 
 by Corollary  \ref{cor2.24}  the  third term  admits this same bound. 
 Thus we see the difference on the left side of (\ref{W-P}) is negligible.  On the other hand  by Proposition \ref{prop3}
$$P_z[\sigma_{[R,\infty)} = \tau_{U(R)}  <\sigma_A]  \frac{2Ry}{\sigma^2n}p^n(y-x) = \frac{g_A^+(x)y}{n}p^n(y-x)\{1+o(1)\}.$$
Consequently   we obtain  (\ref{d_sum}) as required. 
\qed

\v2\v2
{\it Proof of Theorem \ref{thm1}.} \, Owing to  Lemmas \ref{lem4.1} and \ref{lem4.2} we may and do suppose  $|x|\vee |y| \leq M$ for a constant $M$ (cf. Remark 1 (b)).  We make the same argument (with the same $R$) as in the preceding proof
except that in it $y$ is supposed to tend to $+\infty$ and  we have neglected  the contribution of the sum over $z\leq -R$  to the double sum in (\ref{decomp}) but here $y$ remains  bounded and  the contributions  both from $z\geq R$ and from $z\leq -R$  become relevant.    We apply   Lemma \ref{lem4.2} to see that   if $R\leq |z|<\sqrt n$, then  in  the double sum in (\ref{decomp}),
 $$p_A^{n-k}(z,y)= p^{n-k}_{-A}(- y,-z)  = \left\{ \begin{array} {ll}  {\displaystyle \frac{g^-_{-A}(-y)z}{n} p^n(y-z) \{1+ o(1)\} } \quad &\mbox{for}\quad  z\geq R,\\[4mm]
 {\displaystyle \frac{g^+_{-A}(-y)(-z)}{n} p^n(y-z) \{1+ o(1)\} } \quad &\mbox{for}\quad  z\leq -R;
  \end{array} \right. $$
 also by Theorem A (iii), if $|z|\geq \sqrt n$, then $p_A^{n-k}(z,y) \leq C|y|/z^2$. The evaluation of the term  $\e(n,x,y) $ given in (\ref{e})  is valid and  making  estimation of the overshoots as in the last several lines of the preceding proof  we deduce  
\beq
p_A^n(x,y) 
&=& P_x[\tau_{ U(R)} = \sigma_{[R,\infty)}  < \sigma_A]\frac{g^-_{-A}(-y) R}{n} p^n(y-x)\{1+ o(1)\}\\
&&+ \, P_x[\tau_{ U(R)} =\sigma_{(-\infty,- R]} < \sigma_A]\frac{g^+_{-A}(-y) R}{n} p^n(y-x)\{1+ o(1)\} \\
&&  +\, O(n^{-2}). 
\eeq
Substitution from  the formula of Proposition \ref{prop3} and its dual  therefore  concludes the required relation of  Theorem \ref{thm1}.   \qed
\v2
{\sc Remark 5.}\, Corollary \ref{cor1} is verified  virtually in
  the proof above. For in it  we may replace $p_A^{n-k}(z,y)$ by $\widehat p_A^{\,n-k}(z,y)$
with the difference when  $y\in A$ in  which case  the latter  represents the hitting distribution of $A$ in space and time, i.e., $\widehat p_A^{\,n-k}(z,y) = P_z[\sigma_A =n, S_n=y]$.  
A direct derivation   may  also be  made by substituting  the expression of $p_A^{n-1}(x,y)$   given in  Theorem \ref{thm1}   into  the identity
$$P_x[ \sigma_A =n, S_{n} =\xi] =\sum_{y\notin A} p_A^{n-1}(x,y)p(\xi-y),$$
(the last leaving  decomposition) and then using the dual relations of  (\ref{harm_g}) that read
\[
\sum_{y\notin A} p(\xi -y)g^\pm_{-A}(-y) = g^{\pm}_{-A}(-\xi), \quad \xi\in A.
\]

\section{ Refinements in   case $xy<0$}
 In the case $xy<0$ the range of validity of formula (\ref{eq_thm1})  is restricted to  $|x|\wedge |y| <M$ in Theorem \ref{thm1}.  In this section we remove this restriction under  additional conditions on $p$. 
 
 We call
$$D^+_A = \lim_{x\to\infty}[ \sigma^2 a(x)- x - \sigma^2 g^-_A(x)].$$
By (\ref{repr_g-}) the limit exists and for all $\xi_0\in A$,
\beqn\label{D++}
D^+_A =\sum_{\xi \in A} H^{+\infty}_A(\xi)[\sigma^{2} a(\xi -\xi_0) - (\xi- \xi_0)].
\eeqn
$D^+_A$ is positive unless   either $A$ consists of a single point or
the walk $S$ is left-continuous when what we discuss below becomes  trivial. 

\v2
\n
\begin{Prop} \label{prop4.1} \, For any  $M>1$ and $c\in \Z$, 
 as $x\wedge (-y) \wedge n\to\infty$ under the condition  $x\vee(- y) <M\sqrt n$
\beqn\label{4.1}
 p_{\{c\}}^n(x,y)  - p^n_A(x,y) =  \Big(D^+_A +o(1)\Big) \frac{ x+|y|}{\sigma^2 n}p^n(y-x);
\eeqn
  for $c\in A$, in  other words, the probability that the path of a pinned walk  of length   $n$  joining   $x$ and  $y$ passes through $A$ but avoids  the point $c$ is asymptotically equivalent to $D^+_A(x+|y|)/\sigma^2 n$.
\end{Prop}
\v2
\n
\pf\,  
Throughout the proof  the variables $x, y$ and $n$   are  assumed  subject to
the restriction $x\vee |y|<M\sqrt n$ and $y<0<x$.   

 We can suppose  that $c\in A$. Indeed, if  (\ref{4.1}) is valid under this condition, taking  any  $b\notin A $ and applying (\ref{4.1}) with $\{c,b\}$ in place of $A$, we have two equalities and the subtraction of them yields 
 $$p^n_{\{c\}}(x,y) - p^n_{\{b\}}(x,y) = (x+|y|)\times o(n^{-3/2}),$$
  which, combined with the expression for $p_{\{c\}}^n-p_A^n$, 
  shows  (\ref{4.1})  with $b$ in place of $c$. 
   
   Now  suppose $c\in A$.  Then  
\beqn\label{101}
p_{\{c\}}^n(x,y) - p^n_A(x,y)= \sum_{k=1}^n\sum_{\xi\in A\setminus\{c\}} P_x[ \sigma_A=k, S_{k}=\xi] p_{\{c\}}^{n-k}(\xi,y).
\eeqn  
Let $\e$ be a positive  number (small enough that  $\e M^2 <1/2$). For $\xi\in A$, by Corollary \ref{cor1} and (\ref{g_pm}) uniformly  for $k> \e x^2$,  as $x \to \infty$
$$P_x[ \sigma_A=k, S_{k}=\xi] = \frac{\sigma ^2 a(x) +x}{2k}p^{k}(\xi-x)H^{+\infty}_A(\xi)\{1+o(1)\}$$
whereas by Theorem \ref{thm1} (or Theorem A) as $y \to -\infty$, uniformly for $0\leq k\leq  n-\e y^2$,
$$p^{n-k}_{\{c\}}(\xi, y) = \frac{g^-_{\{c\}}(\xi)}{n-k} |y|p^{n-k}(y-\xi)\{1+o(1)\} \quad  \mbox{for}\quad \xi\in A. $$
Denote  the outer  sum in (\ref{101}) restricted on an interval  $a< k\leq b$ by $I_{(a,  b]}$.
We  then  obtain 
\beq
&& I_{( \e x^2, \, n-\e y^2]} \\
&&=\sum_{\xi\in A\setminus \{c\}} H_A^{+\infty}(\xi) g^-_{\{c\}}(\xi) \sum_{ \e x^2 < k\leq n-\e y^2}\frac{xp^k(\xi-x)}{k} \cdot \frac{|y|p^{n-k}(y-\xi)}{n-k}\{1+o(1)\}.
\eeq
Suppose $p^n(y-x) > 0$. Then recalling  that $\nu$ denotes  the period of the walk, we apply the local limit theorem (\ref{llt}) to rewrite the inner sum as
$$\nu 
\int_{\e x^2}^{n-\e y^2} \frac{x{\sf g}_s(x)}{s} \cdot \frac{|y|{\sf g}_{n-s}(y)}{n -s}ds\{1+o(1)\}.$$
Since $ x{\sf g}_s(x)/s$ is the passage time density of Brownian motion, the integral above, if the range of integration is extended to the interval $(0,n)$,  becomes 
$$\frac{x+|y|}{ n}{\sf g}_n(x+|y|) = \frac{x+|y|}{\nu n}p^n(y-x)\{1+o(1)\}.$$
Plainly 
$$\int_0^{\e x^2}x {\sf g}_s(x)s^{-1}ds =\int_0^\e {\sf g}_u(1)u^{-1}du = o(\e)$$
 as $\e\downarrow 0$  and similarly  $\int_{n-\e y^2}^n |y|{\sf g}_{n-s}(y)(t-s)^{-1}ds =o(\e)$, so that for $\e$ small enough,
$$\bigg(\int_0^{\e x^2}+ \int_{n-\e y^2}^n\bigg) \frac{x{\sf g}_s(x)}{s} \cdot \frac{|y|{\sf g}_{n-s}(y)}{n-s}ds \leq \frac{\e}{n}[|y|{\sf g}_n(y) + x{\sf g}_n(x)].
$$
Since   $g^-_{\{c\}}(\xi) = a(\xi-c ) - (\xi-c)/\sigma^2$ for $\xi \neq c$, we also have  $$\sum_{\xi\in A\setminus \{c\}} H_A^{+\infty}(\xi)g^-_{\{c\}}(\xi) = D_A^+/\sigma^2.$$  We
then put these together to  obtain  
\beqn\label{4.2}
I_{ (\e x^2,\, n-\e y^2]} = D_A^+\frac{x+|y|}{ \sigma^2 n}p^n(y-x)\Big[1+o(1) + O(\e)\Big],
\eeqn
where $o(1)\to 0$ as $x\wedge(-y)\wedge n\to \infty$ for each  $\e>0$ and $O(\e)$ is uniform in $x, y$ and  $n$.

As for the remaining parts of the double sum in (\ref{101}) we observe
\begin{eqnarray}\label{eq5.5}
I_{ (0, \,\e x^2]}  &\leq&  C P_x[\sigma_A\leq \e x^2] p_{\{c\}}^n(\xi, y)  \nonumber\\
&\leq& C'P_0[ \max_{k\leq \e x^2} (-S_k)\geq x ]\frac{ |y|}{ n^{3/2}} \leq C''\frac{\e |y|}{n^{3/2}}, 
\end{eqnarray}
where Kolmogorov's inequality is used for the last inequality. 
From  (i) and (iii) of  Theorem A  we derive 
\begin{eqnarray*} 
  p_{\{0\}}^j(\xi,y) &\leq& \frac{C y}{j^{3/2}} e^{-y^2/2\sigma^2 j}  \quad \quad   \quad   \quad  \quad  \quad  \quad \mbox{for} \quad   j>  y^2, \nonumber\\
&\leq&
 \frac{C}{y\sqrt j} e^{-y^2/8\sigma^2 j} + \frac{1 } {y^2}\times o(1) \quad \,\, \mbox{for}  \quad   1\leq j \leq  y^2,
\end{eqnarray*}
where $o(1)$ is bounded and tends to zero as $j\to \infty$ uniformly in $y$,  and  on using these bounds  easy computations yield 
$$I_{ (n-\e y^2,\, n]}  \leq C \bigg[ \sum_{1\leq j \leq \e y^2} \frac{1}{|y|\sqrt j}e^{-y^2/8\sigma^2j}\bigg] \frac{x}{n^{3/2}}<  C' \frac{\e x}{n^{3/2}}.$$
This  together with (\ref{eq5.5}) and  (\ref{4.2}) completes the proof. \qed

\v2
Because of the identity $C_A^+=C^+ - D^+_A$  Theorem \ref{thm2} follows immediately  from Proposition \ref{prop4.1} together with  the result for the case $A=\{0\}$ that is established in \cite[Theorem 1.4]{U1dm}.  (For   verification of the latter one can readily adapt the proof of Proposition \ref{prop4.2} given shortly.) 
  
   Theorem \ref{thm2}  entails  that  (\ref{eq_thm1}), the formula of Theorem \ref{thm1},
  is true uniformly under $|x|\vee |y|= O(\sqrt n)$ and $|x|\wedge |y|= o(\sqrt n)$, if $E[|X|^3; X<0]<\infty$. 
The next  result asserts  that  the same is true  under a much weaker condition on $p$. We state the condition by means of $a$  as follows: 
  for $x, y \geq 1$, 
\beqn\label{H}
\frac{(\sigma^2a(x)-x)/x}{(\sigma^2a(y)-y)/y} \to 0 \quad\mbox{as} \quad \frac{y}{x} \to 0.
\eeqn

\begin{Thm} \label{thm4} \; Suppose  (\ref{H}) to hold true. Then,  for any $M>1$,   uniformly for 
\beqn\label{M_-M}
-M\sqrt n \leq  y<0<  x \le M\sqrt n,
\eeqn
 as $n\to\infty $ and $(x\wedge|y|)/\sqrt n \to 0$  formula (\ref{eq_thm1})  of Theorem \ref{thm1} holds true; in particular if $E[|X|^3; X<0]=\infty$ (so that 
$\sigma^2 a(x)-x \to\infty$ as  $x\to\infty$), as $x\wedge |y|\wedge n\to \infty$ under (\ref{M_-M}) and $x\wedge |y| =o(\sqrt n)$, 
\beqn\label{thm4eq}
p_A^n(x,y)=\frac{(\sigma^2 a(x)-x)|y| + (\sigma^2a(-y) +y) x}{\sigma^2 n} p^n(y-x)\{1+o(1)\}.
\eeqn
\end{Thm}
\v2

If $E[|X|^3; X<0] <\infty$,  then (\ref{H}) holds in view of (\ref{C_+}). If $E[|X|^3; X<0] =\infty$,   by (\ref{Lamb0})
 of Lemma \ref{lem3.1} the condition (\ref{H}) may be  expressed explicitly  in terms of $p$ and  it in particular follows that if $F(x)$ is regularly
varying at $-\infty$
 with an exponent $\a$  less than $-2$, then $\sigma^2  a(x)-x$ varies regularly with  exponent $(\a+3)\vee 0$,  which is less than 1, hence
  (\ref{H})  holds; in the critical case $\a=-2$ (with $\sum F(x)|x|<\infty$) $(\sigma^2a(x)-x)/x$ becomes slowly varying so that (\ref{H}) is violated. 

By virtue  of Proposition \ref{prop4.1}  Theorem \ref{thm4} follows if we show

\v2\n
\begin{Prop} \label{prop4.2} \;  Suppose  (\ref{H}) to hold true. Let $y<0<x$ and suppose  the walk is not left-continuous (so that $p_{\{0\}}^n(x,y)$ is not identically zero).  Then for any $M>1$, uniformly for $x,|y|\le M\sqrt n$, as $x\wedge |y|\wedge n\to \infty$  under $(x\wedge|y|) = o(\sqrt n)$
\beqn\label{prop5_eq}
p_{\{0\}}^n(x,y)=\frac{(\sigma^2 a(x)-x)|y| + (\sigma^2a(-y) +y) x}{\sigma^2 n} p^n(y-x)\{1+o(1)\}.
\eeqn
\end{Prop}
\v2\n
\pf ~ The proof parallels that given for the same formula but under the additional  condition $E[|X|^3; X<0] <\infty$ in \cite{U1dm}. For simplicity we suppose $\nu=1$, i.e.,   the walk is  aperiodic.  Recall that $h_x(n,z)$ is  the space-time  hitting probability of $(-\infty,0]$.  Making decomposition
\beqn\label{eq5.10}
p_{\{0\}}^n(x,y)=\sum_{k=1}^n \sum_{z < 0}h_x(k,z)p_{\{0\}}^{n-k}(z,y),
\eeqn
we break  the double sum into three parts by partitioning the range of the outer summation as follows
\[
1\leq  k< \e n; ~~ \e n \leq k \leq (1-\e)n;~~ (1-\e)n<k\leq n
\]
and call the corresponding sums $I,~II $ and $I\!I\!I$, respectively. Here $\e$ is a positive constant  less than $1/4$ that will be chosen small.

We call
\beqn\label{lamb}
\la(x) := a(x)- \frac{x}{\sigma^2}
\eeqn
and  suppose $\la(x)\to\infty$ as $x\to\infty$, or equivalently  $E[|X|^3; X<0]=\infty$; otherwise  (\ref{prop5_eq}) is already verified in Theorem \ref{thm2}. This supposition permits using 
 (\ref{Hab}), which states
\beqn
\label{Hy}
  \sum_{z=-\infty}^{-1} H^{+\infty}_{(-\infty,0]}(z)( |z|\wedge |y|) 
\sim    \frac{\sigma^2}2 \la(-y)  \quad \quad (y\to-\infty).
\eeqn
Let $x\wedge |y|= o(\sqrt n\,)$.
By duality one may suppose that $y=o(\sqrt n)$. Then the hypothesis (\ref{H}) says
 \beqn\label{hyp}
|y| \la(\sqrt n\,) = \sqrt n\, \la(-y) \times o(1).
 \eeqn
Proposition \ref{prop1.1} implies
$$p_{\{0\}}^{n-k}(z,y)\leq \frac{C(|z|\wedge \sqrt n\,) |y|}{n^{3/2}}\quad\quad   (k\leq \e n, z\leq -1), $$
whereas  using Lemma \ref{thm1.4} [(\ref{h0}, \ref{alpha1})]  as well as (\ref{Hy})   one deduces
\beqn\label{hh}
\sum_{k\geq \e n} \sum_{z<0} h_x(k,z)(|z|\wedge \sqrt n\,) \le M_\e x \frac{\la(\sqrt n\,) }{\sqrt{ n}}. 
\eeqn
  By virtue of  (\ref{hyp})   we therefore  obtain 
$$II \leq M_\e \frac{x|y|}{n^{3/2}}\cdot \frac{\la(\sqrt n\,)}{\sqrt{ n}}= \frac {\la(-y)x}{n^{3/2}}\times o_\e(1)$$
Here (and below) $M_\e$ designates  a constant that may depend on $\e$ but not on the other variables and $ o_\e(1) \to 0$   as $n\to\infty$ and $\e\to 0$ in this order uniformly in $x, y$ (subject to the constraints  in the proposition).

Similarly, on using Theorem A  
\begin{eqnarray}\label{I}
I&=&\sum_{1\le k< \e n}\,\sum_{\,- \sqrt n/\e\,\leq z <0}h_x(k,z)\cdot \frac{\sigma^4a(z)a(-y)+zy}{\sigma^2(n-k)}{\sf g}_{n-k}(y)(1+o_{\e}(1)) \nonumber\\
&& \,+ \sum_{z< -\sqrt n/\e\,} H^x_{(-\infty,0]}(z) \times O\bigg( \frac{y}{n}\bigg).
\end{eqnarray}
To find an upper bound of the last sum we use the identity  (\ref{A.03}), which may be  written as 
\beqn\label{A.02}
\sum_{z<0} H^x_{(-\infty,0]}(z) \la(z)=\la(x).
\eeqn
By   Markov's inequality  and $\la(z)> -z/\sigma^2$ this   entails  
\beqn\label{hyp2}
\sum_{z<- \sqrt n/\e} H^x_{(-\infty,0]}(z) \leq \frac{\e}{\sqrt n} \sum_{z<- \sqrt n /\e} |z|  H^x_{(-\infty,0]}(z) \leq \frac{\e\sigma^2\la(x)}{\sqrt n}.
 \eeqn
For the evaluation of the  double sum in (\ref{I})  we may replace $(n-k)^{-1}{\sf g}_{n-k}(y)$ by $n^{-1}{\sf g}_n(y)(1+O(\e))$.  Since  $y$ is supposed to go to $-\infty$, we may also replace $\sigma^4a(z)a(-y)+zy =\sigma^4\la(z)a(-y)+ \sigma^2z\la(-y) $  by  $\sigma^2\la(z)|y|$ and  in view of (\ref{hh}) we may extend the range of the outer summation    to  the whole
 half line $k\geq 1$. 
Using (\ref{hyp}), (\ref{I}), (\ref{A.02})  and (\ref{hyp2})  we then deduce that  
 $$ I=  \frac{\la(x) |y|{\sf g}_{n}(y)}{n} [1+O(\e) +o_\e(1)].$$

As for $ I\!I\!I$   observe that
 \beqn\label{III0}
  \sum_{(1-\e)n<k\leq n}p_{\{0\}}^{n-k}(z,y)=g(z, y)-r_n(z,y)  \qquad (z<0)
  \eeqn
 with
 $0\leq r_n(z,y) :=  \sum_{j\geq \e n}p_{\{0\}}^j(z,y)  \leq C(|z|\wedge \sqrt n\,)|y|/\sqrt{\e n},$
as deduced from Proposition \ref{prop1.1} for $|z|<\sqrt n$ and by the bound $g(z, y)\leq O(y)$ for  $z|\geq \sqrt n$; hence by (\ref{hh}) and (\ref{hyp})
$$ \sum_{k\geq (1-\e) n} \sum_{z<0} h_x(k,z)r_n(z,y) \leq C\frac{x|y|}{n^{3/2}}\cdot \frac{\la(\sqrt n\,) }{\sqrt{ \e n}} = \frac{x\la(-y)}{n^{3/2}}\times o_\e(1). $$
By (\ref{III0})  it remains to evaluate the sum  $\sum_{z<0} h_x(k,z)g(z,y)$ uniformly over $(1-\e)n<k \leq n$. To this end we apply Lemma  \ref{thm1.4} again. In  formula  (\ref{h0}) of it the error term $\a_n(x,y)$ satisfies
$$\sum_{z<y} \Big[\sup _{(n/2< k\leq n}\a_k(x,z)\Big]g(z,y)  \leq c_M |y|\frac{\la(-\e n)}{\e n}+  \frac{\la(-y)}{y}\times o_\e(1)$$
according to (\ref{alpha2}).
  Now  using ${\sf g}_{k}(x) ={\sf g}_{n}(x)(1+O(\e))$ ($(1-\e)n <k\leq n$) as well as   $y=o(\sqrt n)$   we apply  (\ref{Hy})  and  (\ref{h0})  to  find that
 $$ I\!I\!I= \frac{f^+(x){\sf g}_{n}(x)}{\sigma^2 n}\sum_{z<0}H_{(-\infty, 0]}^{+\infty}(z)g(z, y)[1+o_\e(1)+ O(\e)],$$
hence by the identity  (\ref{A.000})  
  $$ I\!I\!I=  \frac{x\la(-y){\sf g}_{n}(x)}{ n}[1+o_\e(1)+O(\e)].$$
Adding these contributions yields the desired formula, since $\e$ can be made arbitrarily small. ~~~\qed
\v2
  In view of (\ref{Hbb})  $\la(x)/x$ is asymptotically decreasing  (as $x\to \infty$) in the sense that   $\la(x)/x\sim \mu(x)$ with a decreasing $\mu$, as is noted previously. Keeping this in mind one examines   the proof above and deduces the following upper bound  without  assuming (\ref{H}).
 \begin{Prop}\label{prop60} \; There exists a constant  $C$ such that   
\beqn\label{eq_pr7}
p^n_{\{0\}}(x,y) \leq C\frac {x\la(-y) + |y|\la(x) }{n^{3/2}}\quad\quad (-M\sqrt n <y<0< x<M\sqrt n).
\eeqn
\end{Prop} 
\v2\n

\v2
{\sc Remark 6.} One can  show that $p^n_{\{0\}}(x,y)$ may be of smaller order of the right side of (\ref{eq_pr7}) for suitably chosen  $x, y$ with $x\vee |y| =o(\sqrt n\,)$ if e.g. $P[X<-z] \sim 1/z^2(\log z)^2$ ($z\to +\infty$). Thus formula (\ref{eq_thm1}) is not  generally true for $|x|\vee|y|<M\sqrt n$,  $xy<0$.
\v2
The next result provides a lower bound without  assuming (\ref{H}).

\begin{Prop}\label{prop6}  Given  $M\geq 1$,  for  $-M\sqrt n <y<0< x<M\sqrt n$ satisfying  (\ref{nat_c}), 
$$p^n_{A}(x,y) \geq c\bigg( \sum_{w=2}^{x\wedge|y|} p(-w)w^3\bigg)\frac{x + |y|}{\sigma^4 n}p^{n}(y-x),$$
where $c$ is a positive constant depending only on $M/\sigma$.
\end{Prop}
\v2
If $F(y)$ is regularly varying as $y\to -\infty$ with exponent $ \a\in [-3,-2]$ and  $E[|X|^3;X<0] =\infty$, then in the non-critical case  $\a\neq -2$ we have $\sum_{w=2}^{x} p(-w)w^3 \asymp \la(x)$ and hence  owing to (\ref{H})
$$(x + |y|)\sum_{w=2}^{x\wedge |y|} p(-w)w^3\,  \asymp \, x\la(-y)+ |y|\la(x) \quad\mbox{as}\quad x\wedge (-y) \to\infty$$
so that the above lower bound is exact,
 whereas in the critical case $\a=-2$ when  $w\int_w^\infty F(-t)dt$ is  slowly varying at infinity  
    $\sum_{w=2}^{x\wedge (-y)} p(-w)w^3$ $ =o(\la(x))$ ($x\to\infty$)---as is inferred  by using (\ref{Lamb0})---so that the lower bound is consistent to what is mentioned in Remark 6.
\v2
{\it Proof of Proposition \ref{prop6}.}  In view of Proposition \ref{prop4.1} as well as Theorems \ref{thm2} we may suppose $A=\{0\}$. We also suppose $\nu=1$ for simplicity and  $x\leq |y|$ for reason of  duality. Let  $c_1, c_1', c_2$ etc. denote  positive constants depending only on $M/\sigma$.  Substitution from Theorem C and Theorem \ref{thm1} into (\ref{eq5.10}) shows that for a  constant  $\de>0$,
$$p^n_{\{0\}}(x,y) \geq c_1 \sum_{\de x^2\leq k \leq n/2}\; \sum_{1\leq w \leq \sqrt{k/\de}}\;\sum_{z=-x}^{-1}\frac{xw}{\sigma^2k}p^k(w-x)p(z-w)\frac{zy}{\sigma^2(n-k)}p^{n-k}(y-z).
$$
We take
$$\de = (4 M^2)^{-1}$$
so that $n/2(x-y)x> \de$ hence  $\de x^2 =nx/2(x-y) \leq n/4$ due to   $x\leq |y|$.
 For $k, w, z$ taken from the range of summation above, in view of the local limit theorem
\beq
p^k(w-x)p^{n-k}(y-z) &\geq& c_2 {\sf g}_{k}(w-x){\sf g}_{n-k}(y-z)\\
&\geq& c_2 {\sf g}_k(x){\sf g}_{n-k}(y).
\eeq
 If $0<x <M\sqrt n$ and  $k\geq \delta x^2$,  then  $x\leq \sqrt{k/\delta}$ so that   ${\sf g}_k(x)\leq c_2' e^{-1/2\de\sigma^2}k^{-1/2}/\sigma$; and  also 
  ${\sf g}_{n-k}(y) \leq  {\sf g}_{n/2}(y)\leq c_2'' p^n(y-x)$ ($k\leq n/2$). Hence,    putting
$$m(x) =  \sum_{w=1}^x\;\sum_{z=1}^x p(-z-w)wz
$$
we have
$$p^n_{\{0\}}(x,y) \geq c_3 \frac{m(x)p^n(y-x)x|y|}{\sigma^4  n}\sum_{\de x^2\leq k \leq n/2}\;\frac1{\sigma k^{3/2}}.$$
Since $\de x^2 \leq n/4$, the last sum is bounded below by  $(2-\sqrt 2)/x\sigma \sqrt{\de}$. 
   Finally an easy computation yields  that $m(x)\geq \sum_{w=1}^x\;\sum_{z=1}^{x-w} p(-z-w)wz
\sim  \frac16  \sum_{w=2}^x p(-w)w^3$ ($x\to\infty$), which concludes the proof. \qed

\begin{Cor}\label{cor8}  For each  $M\geq 1$, it holds  under the constraint   $-M\sqrt n <y<0< x< M\sqrt n$  that 
\beq
&& P_x[S_{\sigma_{(-\infty,0]}}<-\eta \,|\,\sigma_{\{0\}}>n, S_n =y] 
\\[2mm]
&&\quad  \longrightarrow  \left\{ 
\begin{array}{ll}
0 \quad\mbox{ as}  \;\;  \eta \to\infty\;\;\;\quad  \mbox{uniformly for} \;\; x, \,y\; &\mbox{ if}\quad E[|X|^3; X<0] <\infty, \\[1mm]
1 \quad\mbox{ as}  \; \;  x\wedge (-y) \to\infty \;\; \mbox{ for each}\; \eta>0\;\; &\mbox{ if}\quad E[|X|^3; X<0] =\infty.
\end{array} \right.
\eeq

\end{Cor}
\v2\n
\pf\,   If  $E[|X|; X<0] <\infty $, then by  $H_{(-\infty,o]}^x(z)\leq CH_{(-\infty,0]}^\infty(z)$ it follows that  for any $w<0$ 
$$\sum_{z<-\eta} H_{(-\infty,-\eta]}^x(z)(|z|\wedge |w|)\leq C\sum_{z<-\eta} H_{(-\infty,0]}^{+\infty}(z)|z| \to 0 \quad (\eta\to\infty),$$ 
and following the proof of Proposition \ref{prop4.2} one can  readily  show the first relation. The second one  follows from Proposition \ref{prop6} since the contribution to the sum (6.10) from   $-\eta \leq z<0$ is $O(\eta (x\vee |y|)/n^{3/2})$ as is easily verified (use $\sum_{k\geq n/2} p_{\{0\}}^{n-k}(z,y)  \leq g_{\{0\}}(z,y) \leq C\eta$ for $z \leq -\eta$), hence  negligible relative to the lower bound given by Proposition \ref{prop6}. \qed

\section{Appendices}

{\bf A.}\;  {\sc A consequence of duality.} 
\v2

In (\ref{greenft}) we have brought in  Green's function 
\beqn\label{g/a}
g(x,y) = a(x)+ a(-y)- a(x-y)
\eeqn
on the space $\Z\setminus \{0\}$. 
By means of $g$  the identity (\ref{2_2}) may be written as
\begin{eqnarray}\label{2_21}
g_A(x,y) &=&   a^\dagger(x-\xi_0) - a(x-y) + E_x[a(S_{\sigma(A)}-y)] - E_x[a(S_{\sigma(A)}-\xi_0)] \nonumber\\
&=&  \de_{x,\xi_0} +g(x-\xi_0, y-\xi_0) - E_x[g(S_{\sigma(A)}-\xi_0, y-\xi_0)]
\end{eqnarray}
(for any  $\xi_0\in A$). It is noted  that  for all $x,y$ 
\[
g_A(x,y) = g(x,y) +O(1).
\]

  We consider the dual  walk, denoted by $\widehat S_n$, that  is a  random walk with transition probability
 $$P[ \widehat S_n =y\,|\, \widehat S_0=x] = p^n(x-y).$$
The objects associated with this walk  are denoted by  $\widehat P_x, \widehat E_x, \widehat p^{\,n}_A$, etc.
Then for $x, y \notin A$, 
\beqn\label{dual_p}
\widehat p^{\,n}_A(y,x) = p^n_A(x,y) = p^n_{-A}(-y,-x),
\eeqn
 the law of $(\widehat S_n)$ with $\widehat S_0=y$  being the same as that of the walk $(-S_n)$:
 $$\widehat p^{\,n}_A(y,x)=P[ -S_1\notin A,\ldots,-S_n \notin A,-S_n = x \,| -S_0=y] = p^n_{-A}(-y,-x). $$
 Hence $ \widehat g_A(y,x) = g_{A}(x,y) =g_{-A}(-y,-x)\; (x, y \notin A)$.
It is noted that the second equality in (\ref{dual_p})  follows directly from  the invariance of the law of the walk under a reversal of the order of the increments $X_1,\ldots, X_n$.

We prove the following result of which the first relation entails   (\ref{C-}).
  \begin{lem} \label{lem7.1}\;  $ \lim_{x\to \infty} [g^-_A(x)- g^-_{-A}(x)] =0$ and  $\,  \lim_{y\to -\infty} [g^+_A(y)- g^+_{-A}(y)] =0$.
\end{lem}
   \v2 \n
 \pf\,   It suffices to show the first relation, the second being its dual. We may suppose $0\in A$.  Using the identity $g_{-A}(x,y) =g_{A}(-y,-x)$,  (\ref{2_1}) and
 (\ref{2_2}) in turn we see that
 \beq
 g_A(x,y) - g_{-A}(x,y) &=& u_A(x)-u_A(-y) + E_x[a(S_{\sigma(A)} -y)] -E_{-y}[a(S_{\sigma(A)} +x)]\\\
 &=& a(x)- E_x[a(S_{\sigma(A)})] -a(-y) +E_{-y}[(S_{\sigma(A)})] \\
 && +\; E_x[a(S_{\sigma(A)} -y)] -E_{-y}[a(S_{\sigma(A)} +x)],
\eeq
Letting $y\to   -\infty$ we obtain
$$g_A^-(x) - g_{-A}^-(x)  =  a(x) - E_x[a(S_{\sigma(A)})] + \sum_\xi H_A^{+\infty}[a(\xi) + \xi/\sigma^2-a(\xi+x)],$$
which tends to zero as  $x\to +\infty$ as desired. \qed
 \v2\v2\n
{\bf B.}\; {\sc Some inequalities concerning $a(x)-x/\sigma^2$.}
\v2

Let  $\la(x)$ be the function defined by (\ref{lamb}) and $\widehat \la(x)$  its dual:
$$\la(x)= a(x) -\frac{x}{\sigma^2} \quad \mbox{and}\quad \widehat\la(x) = a(-x) -\frac{x}{\sigma^2}.$$
Here we collect several formulae satisfied by  $\la(x)$ or $\widehat \la(x)$, all of which rest on  the known results presented in  Section 2; they are  trivial when  $a(x) -|x|/\sigma^2$ converges to a finite number as $x\to\infty$ or $x\to -\infty$ according to the situation.  
It is noted that
$$\la(x+y) -  \la(y) =    a(x+y) -a( y) - {x}/{\sigma^2}$$
and 
$$ \widehat\la(x+y) - \widehat\la(y) =    a(-x-y) -a(-y) - {x}/{\sigma^2}.$$

\v2
(i)\,   {\it For $x>0$ and $y\geq 0$, 
 \beqn\label{A.1}
- \la(x)\times o(1) \leq   \la(x+y) -\la(y)  \leq  \la(x)
  \eeqn 
 and  
 \beqn\label{A.2}
 - \widehat \la(x) \times o(1)   \leq   \widehat\la(x+y) -\widehat\la(y)  \leq    \widehat\la(x), 
 \eeqn
  where $o(1)$ is bounded and, as $y\to\infty$,  tends to zero   uniformly in $x$. }
\v2
The relation (\ref{A.2}) is a dual of (\ref{A.1}).  For the proof of (\ref{A.1}) we may suppose that the walk is not left continuous so that $\la(x) >0$ for $x>0$.
 The second inequality of (\ref{A.1})  is the same as $g(x,-y)\geq 0$. Using  identity    (\ref{A.03}) as well as $g(x,-y)= \sum H^x_{(-\infty,0]}(z) g(z,-y)$ we observe
 \begin{eqnarray}\label{3eqty}
 -\la(x+y)+\la(y) &=&  g(x,-y) - [a(x)- x/\sigma^2]\nonumber\\
&=& \sum H^x_{(-\infty,0]}(z) \Big(g(z,-y) - [a(z)- z/\sigma^2]\Big)\nonumber\\
&=& \sum H^x_{(-\infty,0]}(z) [a(y) -a(y+z) + z/\sigma^2] \nonumber\\
&\leq& \sum H^x_{(-\infty,0]}(z) [a(-z) + z/\sigma^2].
 \end{eqnarray}
The left-most member  tends to zero as  $y\to\infty$ for each $x$  while  the right most is positive
because of the non-left-continuity assumption. Hence   the first inequality of (\ref{A.1}) holds true  for each $x$.  To see the uniformity in $x$ observe that    $a(-z) + z/\sigma^2=o(z)$   as $z\to-\infty$ whereas  $a(z) - z/\sigma^2\geq |z|/\sigma^2$,  then you may apply  
  (\ref{A.03}) again,   in the case when $\sum H^\infty_{(-\infty,0]}(z) |z|=\infty$.  In the other case we use the expression next to the last  in  (\ref{3eqty}) and apply the dominated convergence (or use the fact
  mentioned around (\ref{C_+}) from the outset). 

  \v2
 
 (ii) \, {\it For $0<x\leq y$, 
 \beqn\label{A.3}
   - \widehat\la(x) \leq   \widehat\la(y-x) - \widehat\la(y)  \leq  \widehat\la(x)\times o(1),
\eeqn
where $o(1)$ is bounded and, as  $y-x\to\infty$, tends to zero   uniformly in $x$. }
\v2
On putting $y'=y-x$ and writing 
 $\widehat\la(y-x) - \widehat\la(y)  = \widehat\la(y') - \widehat\la(y'+x),$  there follows from (\ref{A.2}) the formula (\ref{A.3}) with the second inequality restricted to the case when 
  $y' \to \infty$. 
\v2

(iii)\, {\it If $ |x|\leq y$, }
\beqn\label{A.4}
 | \,\widehat\la(y-x) - \widehat\la(y)| \leq C \widehat\la(|x|).
\eeqn 
\v2
 This follows from (\ref{A.2}) if  $x<0$ and from  (\ref{A.3}) if  $x\geq 0$.

\v2\n
{\bf C.}\;  {\sc Comparison between  $p^n_{A}$ and $p^n_{\{0\}}$.}
 
\v2
 What is asserted in Remark 2 for $x, y$ subject  to $|x|\wedge |y| =O(1)$,  not immediate from Theorem \ref{thm1}, is shown in this appendix. As in Remark 2 let $0\in A$ and  $\sharp A\geq 2$ (one may suppose $\sharp A=2$ if he wishes) and    suppose  the walk is aperiodic  and  $|x|\vee |y|=O(\sqrt n)$. 
 
 Suppose that  $y\notin A$ and  $p^n_{\{0\}}(x,y) >0$ for some $n\geq 1$  (hence for all sufficiently large $n$). Then 
 for each (admissible) $x$ and $y$,   the ratio $p_A^n(x,y)/ p_{\{0\}}^n(x,y)$  converges to 
\beqn\label{A.C1}
\frac{g^+_A(x)g^-_{-A}(-y)+ g^-_A(x)g^+_{-A}(-y)}
{g^+_{\{0\}}(x)g^-_{\{0\}}(-y)+ g^-_{\{0\}}(x)g^+_{\{0\}}(-y)}
\eeqn
 as $n\to\infty$, and this ratio
converges to 
\beqn\label{A.C2}
\frac{g^-_{-A}(-y)}{g^-_{\{0\}}(-y)}
 =1-\frac{E_{-y}[\sigma^2a(S_{\sigma(-A)})- S_{\sigma(-A)}]}{\sigma^2a(-y)+y}
 \eeqn
 as $x\to +\infty$ for each $y$ (with $\sigma^2a(-y)+y\neq 0$)  and 
\beqn\label{A.C3}
\frac{ g^-_A(x)}{ g^-_{\{0\}}(x)}
=  1-\frac{E_{x}[\sigma^2 a(S_{\sigma(A)})+ S_{\sigma(A)}]}{\sigma^2 a^\dagger(x)+x} \quad
\eeqn 
as  $y\to -\infty$ for each $x$ (with $\sigma^2a^\dagger(x)+x\neq 0$) and similarly for the cases $x\to-\infty$ and $y\to +\infty$. What is presently  required   is to prove that the ratios on the left sides in (\ref{A.C1}) and in (\ref{A.C2}) are less than unity unless $p^n_A(x,y) = p_{\{0\}}^n(x,y)$ for all $n$.
 It suffices to show that those in (\ref{A.C2})  and  (\ref{A.C3}) are less than unity,  this entailing that  the ratio (\ref{A.C1}) is  also less than unity since $ g^\pm_A(x)\leq  g^\pm_{\{0\}}(x)$,.

For reason of duality we suppose $x\geq 0$ and $|y|<M$ for an arbitrarily chosen  $M>1$. 
We have $p_A^n(x,y) < p_{\{0\}}^n(x,y)$ for some $n\geq 1$  if and only if the walk from $x$ to $y$ can pass through $A$ before visiting $0$ with a positive probability, provided $xy\neq 0$. For $x>0$, this is the case if and only if    the following condition holds
$$
(*) \; \left\{ \begin{array} {ll}  y>0 \;\mbox{ and} \;\; A\cap [1,\infty) \neq \emptyset  \qquad  &\mbox{ if \, $S$ \,is  left-continuous,}   \\[3mm]
\left. \begin{array} {ll}\!\!\! \mbox{either}  \;\;\; y>0 \;\mbox{ and} \; A\cap [1,\infty) \neq \emptyset \\ 
\!\!\! \mbox{or}\qquad \;  y<0 \;\;\mbox{and $S$ is not left-continuous}\end{array} \right\} \quad &\mbox{ if  \, $S$ \,is  right-continuous.}
 \end{array}\right.  $$
\v2\n 
 The case $x=0$ may be reduced to the case $x\neq 0$ and obviously $p_A^n(0,y) < p_{\{0\}}^n(0,y)$ ($y\neq 0$) for some $n\geq 1$ under $(*)$. Now, using relation (\ref{a_rl}), one can  easily  check that the required inequality is true.

\v2\v2
{\bf Acknowledgments.}\, 
 I wish to thank  anonymous  referees for their   helpful comments that motivate the author to
 improve the paper in various aspects in  revision.


\begin{thebibliography}{99}
\baselineskip=9pt


\bibitem{BD}  A. Bryn-Jones  and   R. Doney,  A functional limit theorem for random walk conditioned to stay positive, J. London Math. Soc. (2) {\bf 74} (2006), 244-258


\bibitem {Car} F. Caravenna,  A local limit theorem for random walks conditioned to stay positive, Probab. Theory Rel. Fields {\bf 133} (2005), 508-530



\bibitem {CMM}  P. Collet, S. Martinez and J. Martin,  Asymptotic behaviour of a Brownian motion on exterior domains, Probab. Theo.  Rel. F. {\bf 116} (2000), 303-116.

\bibitem {D} R.A. Doney, Local behaviour of first passage probabilities, Probab. Theor. Rel. Fields, {\bf 152}, (2012), 559-588. 



\bibitem {H}  G. Hunt, Some theorems concerning Brownian motion, Trans. Amer. Math. Soc. {\bf s1}, (1956),  294-319.



\bibitem {KSK} J.G. Kemeny, J.L. Snell and A.W. Knapp, Denumerable markov chains, D. Van Nostrand, Co, Inc., Princeton, N.J., 1966.


\bibitem {K} H. Kesten,  Ratio limit theorems II, Journal d'Analyse Math. {\bf 11}, (1963), 323-379.

\bibitem {LL} G.F. Lawler,  V. Limic, Random walk:  a modern introduction, Cambridge univ. press, Cambridge, 2010.  


\bibitem{PS} S. C. Port and  C. J. Stone, Hitting time and hitting places for non-lattice recurrent 
random walks, Jour. Math. Mech. {\bf 17},  (1967), 35-57. 

\bibitem {S} F. Spitzer,  Principles of Random Walk, 2nd ed.,  Springer, New York, 1976. 



\bibitem{Ufh}  K. Uchiyama, The first  hitting time of a single point for random walks, Elect. J. Probab. {\bf 16},  (2011), 
1160-2000.

\bibitem{U1dm}  K. Uchiyama, One dimensional lattice random walks with absorption at a point/on 
a half line, J. Math. Soc. Japan {\bf 63}, (2011),  675--713. 

\bibitem{Uf.s}  K. Uchiyama,  Asymptotic behaviour of  a random walk  killed on a finite set, Potential Anal.  {\bf 44}, (2016), 497-541. 

\bibitem{Ubdsty} K. Uchiyama, The   transition density  of Brownian  motion  killed on  a bounded set, preprint, available at: http://arxiv.org/abs/1603.03902

\bibitem{U1dmk}  K. Uchiyama, One dimensional  random walks killed on a finite set II, preprint

\bibitem{VW}  V. A. Vatutin and V. Wachtel, Local probabilities for random walks conditioned to stay positive, Probab. Theory Rel. Fields, {\bf 143} (2009), 177-217

\bibitem{Ws}  G. Woess, Random walks on infinite graphs and groups, Cambridge tracts in mathematics 138, Cambridge Univ. Press,  Cambridge, 2000.

\end{thebibliography}
\end{document}